\documentclass[10pt]{article}
\usepackage{amstext,amsthm,amssymb,amsmath}
\usepackage{graphicx}
\usepackage{subfigure}

\usepackage{color}

\usepackage{fullpage}
\usepackage{multirow}
\usepackage{tabulary}

\hfuzz=\maxdimen
\title{Hierarchical multiscale modeling for flows in fractured media using Generalized Multiscale Finite Element Method}

\author{Yalchin Efendiev \thanks{Department of Mathematics \& Institute for Scientific Computation (ISC),
Texas A\&M University,
College Station, Texas, USA
and Center for Numerical Porous Media (NumPor),
King Abdullah University of Science and Technology (KAUST),
Thuwal 23955-6900, Kingdom of Saudi Arabia. Email: {\tt efendiev@math.tamu.edu}.}
\and
Seong Lee\thanks{Chevron ETC, Houston, TX 77002}
\and
Guanglian Li\thanks{Department of Mathematics, Texas A\&M University, College Station, TX 77843-3368}
\and
Jun Yao\thanks{School of Petroleum Engineering, China University of Petroleum (East China), Qingdao City, Shandong Province, 266555, China}
\and
Na Zhang\thanks{School of Petroleum Engineering, China University of Petroleum (East China), Qingdao City, Shandong Province, 266555, China}
}

\begin{document}
\maketitle

\begin{abstract}
In this paper, we develop a multiscale finite element method for solving flows
 in fractured
media. Our approach is based on Generalized Multiscale Finite Element Method (GMsFEM),
where we
represent the fracture effects on a coarse grid via multiscale basis functions.
These multiscale basis functions are constructed in the offline stage via local
spectral problems following GMsFEM. To represent the fractures on the fine grid, we consider
two approaches (1) Discrete Fracture Model (DFM) (2) Embedded Fracture Model (EFM)
and their combination.
In DFM, the fractures are resolved via the fine grid, while in EFM the fracture and
the fine grid
block interaction is represented as a source term. In the proposed multiscale method,
additional multiscale basis functions are used to represent the long fractures, while
short-size fractures are collectively represented by a single basis functions.
The procedure is automatically done via local spectral problems. In this regard, our approach
shares common concepts with several approaches proposed in the literature as we discuss.
Numerical results are presented where we demonstrate how one can adaptively add basis
functions in the regions of interest based on error indicators. We also discuss the use
of randomized snapshots (\cite{randomized2014}) which reduces the offline computational cost.

\end{abstract}

\section{Introduction}

Many porous media flow and transport processes are dominated by the presence
of fractures. The fractures present high conductivity conduits
and their effects need to be captured accurately in the simulations.
Because the fractures have very small width,
they are often represented as zero thickness objects
in numerical simulations. Because there are many
fractures of different lengths (connected or not connected),
the problem is inherently multiscale and
efficient multiscale techniques are needed
to represent their effects.

In this paper,
we present an application of Generalized Multiscale
Finite Element Method for flows in fractured media.
The study of flow in fractured media on the fine grid has
been conducted in numerous papers.
These models
applied various fracture models, such as the Discrete Fracture Model (DFM),
Embedded Fracture Model (EFM), the single-permeability model, and
the multiple-permeability models (\cite {wu2011multiple,Baca84, Karimi-Fard03, Lee01,Hajibeygi11(2),shu06}).
Though these approaches are designed for fine-scale simulations,
a number of these approaches represent the fractures at a macroscopic level.
For example, multiple-permeability models represent the network of connected
fractures macroscopically by introducing several permeabilities in each block.
The EFM (\cite{Lough97, Li08, Lee01}) models the interaction of fractures with
the fine-grid blocks separately for each block.
These approaches can be generalized by incorporating the interaction
of fractures and permeability heterogeneities locally
which can lead to efficient upscaling techniques. Some general upscaling
techniques for flows in fractured network
are presented in \cite{Durlofsky91, gong08},
where
the authors introduce the fracture and matrix interaction parameters
hierarchically.



In recent papers \cite{hajibeygi2011loosely},
several multiscale approaches
are proposed for representing the fracture effects. These approaches share
common concepts with the methods that we discuss here in a sense that they
add new degrees of freedom to represent the fractures on a coarse grid.
The main difference is that our approaches use local spectral problems
accompanied by adaptivity to detect the region where to add new basis functions.
In this regard, the procedure of finding multiscale basis functions
and the enrichment procedure is different from existing techniques.

In this paper, we use a general multiscale finite element framework, GMsFEM.
GMsFEM is a flexible general framework that generalizes the Multiscale
Finite Element Method (MsFEM) (\cite{hw97})
by systematically enriching the coarse spaces and taking into account small scale information and complex input spaces.
In this work, we use Discrete Fracture Model (DFM)
(\cite{Noorishad82, Hoteit08, Zhang12, Huang11}) for the simulation
of the fine-scale problem as well as in the construction of spectral
problem for the GMsFEM.
In GMsFEM approach, as in many
multiscale model reduction techniques, which divides the computation into
two stages: the offline stage and the online stage. In the offline stage,
a small dimensional space is constructed that can be
used in the online stage to construct multiscale basis functions.
These multiscale basis functions can be re-used for any input parameter
to solve the problem on a coarse grid. The main idea behind the construction
of offline and online spaces is the selection of local spectral problems
and the selection of the snapshot space.
In \cite{egh12}, we propose several general strategies. The main idea of
this paper is to combine GMsFEM
with DFM and/or EFM to solve the flow problem in fractured media.
We present a multiscale basis construction, adaptivity,
and the use of randomized snapshots
to reduce the computational cost.


Our approaches share some common concepts with hierarchical fracture modeling
proposed by Lee et al. \cite{Lee01}, where the main idea is to
homogenize small-length fractures
(with the length smaller than the coarse block),
while represent the large-length fractures.
Our methods via local multiscale basis functions homogenizes
small-size fractures (by lumping their effect into a single-per-node multiscale basis
function) and
represent the long-size fractures in each coarse block.
A global coupling, such as finite
element in this case, recover an accurate representation of long size fractures by
coupling all the information together.


We present several numerical examples to illustrate the performance of the proposed approach.
We consider the fine-scale fracture representations using DFM, EFM as well as
coupled DFM and EFM, where the shorter fractures are represented by DFM and the longer ones
EFM. All the numerical results indicate that the proposed GMsFEM
is robust and accurate. Besides, we test the performance of
adaptive algorithm in \cite{Chung2013} for the problem in this paper.
To reduce the computational cost, we propose the use of randomized
snapshots where only a few snapshots are constructed and used to
construct multiscale basis functions. We would like to emphasize that
our goal is to develop and show the performance of GMsFEM for
flows in fractured media and we do not make
any comparisons between different fine-scale fracture models.

The rest of the paper is organized as follows. In Section \ref{prelim}, we present
preliminaries. The fine-scale fracture modeling techniques are
briefly reviewed in Section \ref{sec:DFN}. The construction of the coarse
spaces for the GMsFEM is displayed in Section \ref{sec:GMsFEM}.
In Section \ref{sec:numerical},
numerical results for several representative examples are showed.
Finally, we conclude our paper with some remarks
in Section \ref{sec:conclusion}.

\section{Preliminaries}
\label{prelim}
\begin{figure}[htb]
  \centering
  \includegraphics[width=0.65 \textwidth]{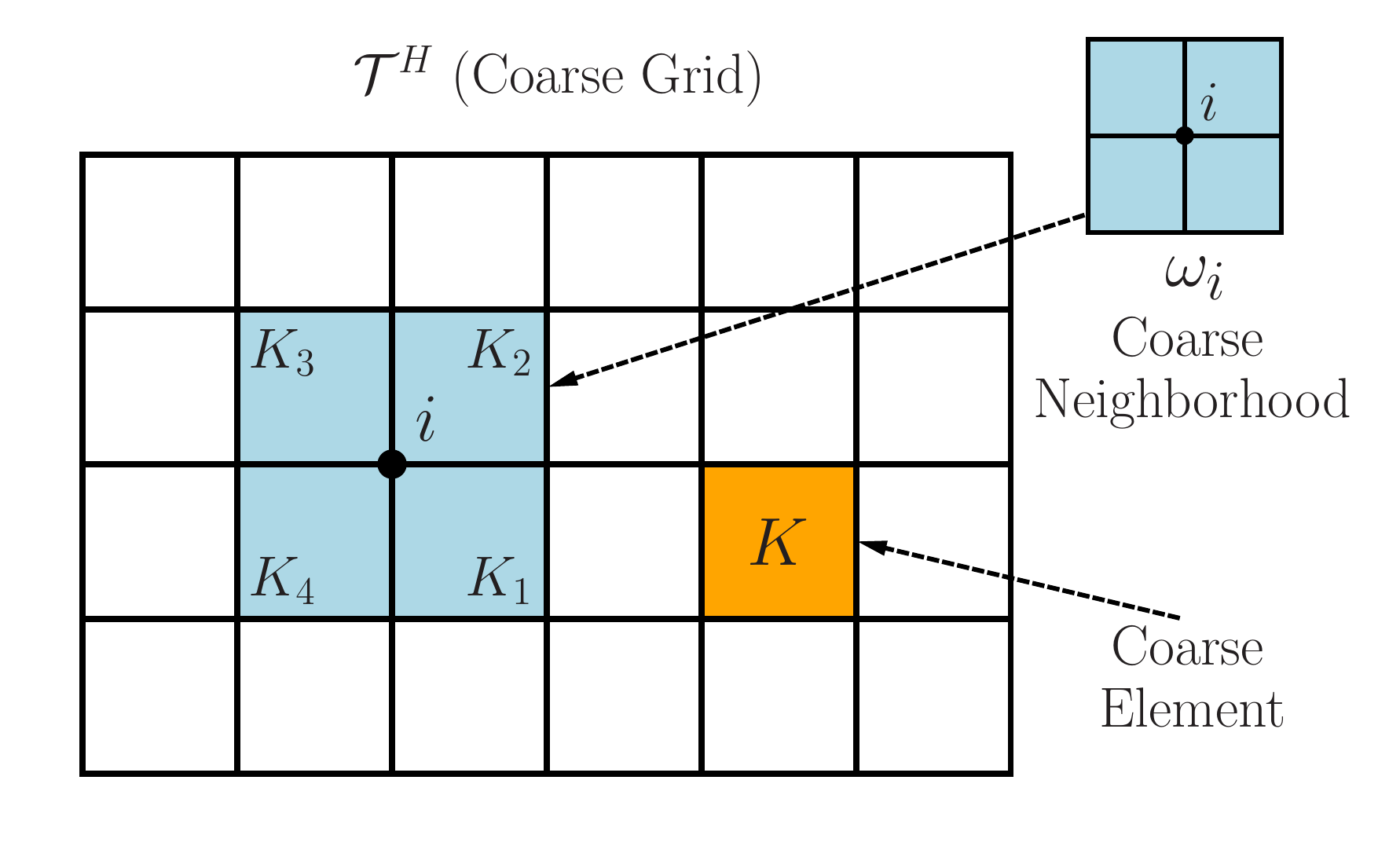}
  \caption{Illustration of a coarse neighborhood and coarse element}
  \label{schematic}
\end{figure}

\begin{figure}[htb]
  \centering
  \includegraphics[width=0.40 \textwidth]{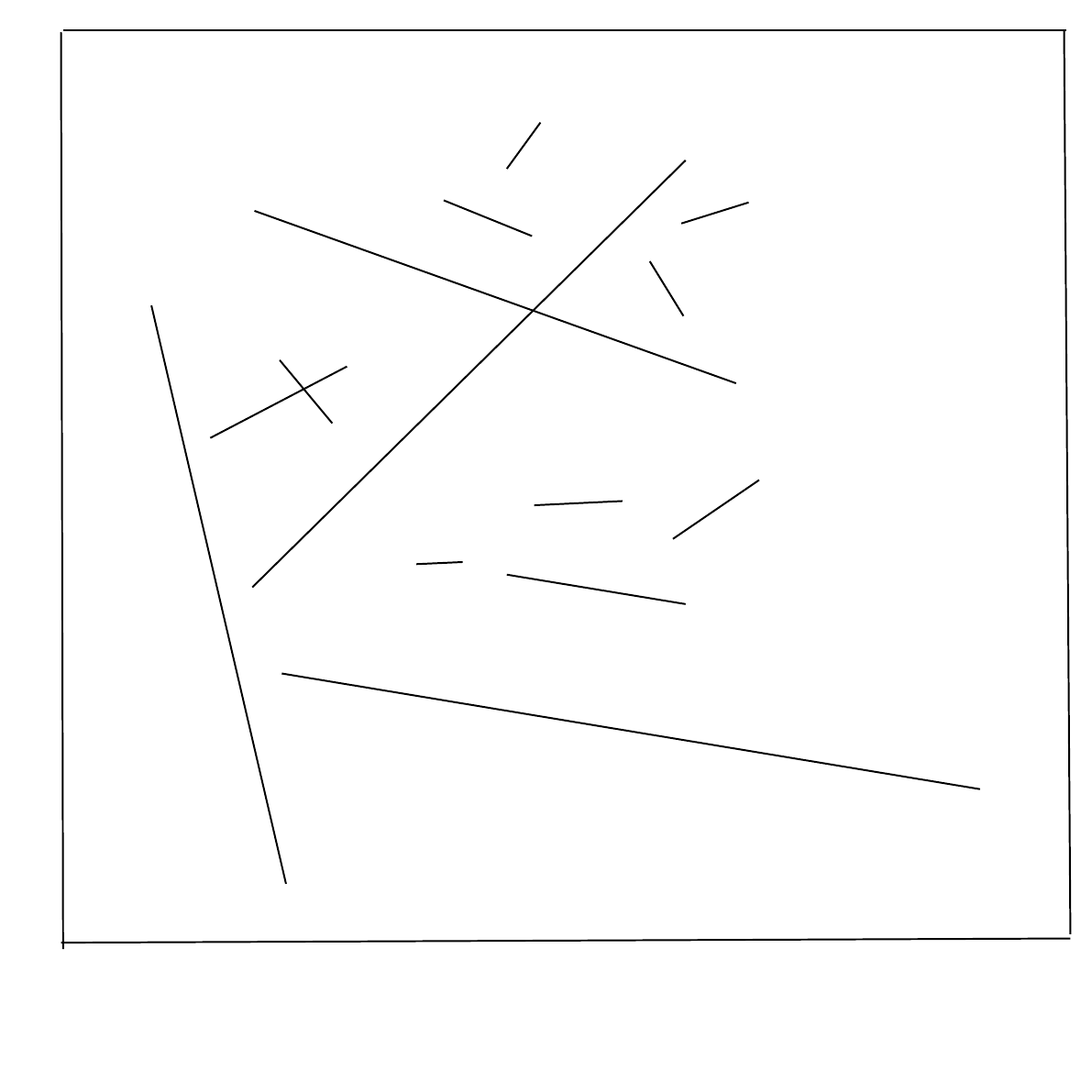}
  \caption{A heterogeneous fracture field.}
  \label{fig:fracture_field}
\end{figure}

In this paper, we study high-contrast flow problem
\begin{equation} \label{eq:original}
-\mbox{div} \big( \kappa(x) \, \nabla u  \big)=f \quad \text{in} \quad D,
\end{equation}
where $\kappa(x)$ has fractures with high values and low values in the matrix
(see Figure \ref{fig:fracture_field} for illustration).
To discretize (\ref{eq:original}), we introduce the notions of fine and coarse grids.
Let $\mathcal{T}^H$ be a usual conforming partition of the computational domain
$D$ into finite elements (triangles, quadrilaterals, tetrahedra, etc.). We denote this
partition as the coarse grid and assume that each coarse element is partitioned into
 a connected union of fine-grid blocks. The fine grid partition will be denoted by
$\mathcal{T}^h$, and is by definition a refinement of the coarse grid $\mathcal{T}^H$.
We use $\{x_i\}_{i=1}^{N}$ (where $N$ denotes the number of coarse nodes) to denote the vertices of
the coarse mesh $\mathcal{T}^H$, and define the neighborhood of the node $x_i$ by
\begin{equation} \label{neighborhood}
\omega_i=\bigcup\{ K_j\in\mathcal{T}^H; ~~~ x_i\in \overline{K}_j\}.
\end{equation}
See Figure~\ref{schematic} for an illustration of neighborhoods and elements subordinated to the coarse discretization. We emphasize the use of $\omega_i$ to denote a coarse neighborhood, and $K$ to denote a coarse element throughout the paper.

Next, we briefly outline the GMsFEM. We will consider the
continuous Galerkin (CG) formulation and signify $\omega_i$ as the support of basis functions.
We denote the basis functions by $\psi_k^{\omega_i}$, which is supported in $\omega_i$,
and the index $k$ represents the numbering of these basis functions.
In turn, the CG solution will be sought as $u_{\text{ms}}(x)=\sum_{i,k} c_{k}^i \psi_{k}^{\omega_i}(x)$.
Once the basis functions are identified (see the next section),
the CG global coupling is given through the variational form
\begin{equation}
\label{eq:globalG} a(u_{\text{ms}},v)=(f,v), \quad \text{for all} \, \, v\in
V_{\text{off}},
\end{equation}
where  $V_{\text{off}}$ is used to denote the space spanned by those basis functions
and $a(\cdot,\cdot)$, $(f,\cdot)$ are bilinear (or linear) form corresponding to (\ref{eq:original}) defined by
$ \displaystyle a(u, v) = \int_D \kappa(x) \nabla u \cdot \nabla v \, dx$, and $ \displaystyle (f,v) = \int_D f v \, dx$.
We also note that one can use
discontinuous Galerkin formulation (see e.g., \cite{efendiev4(1),WaveGMsFEM,eglmsMSDG}) to couple multiscale basis functions
defined on $K$.

Let $V$ be the conforming finite element space
with respect to the fine-scale partition $\mathcal{T}^h$.
We assume $u\in V$ is the fine-scale solution satisfying
\begin{align}\label{eqn:fine-scale prb}
a(u,v) = (f,v), \quad v\in V.
\end{align}
In the next section, we will introduce fine-grid discretizations for fractures.
As we have mentioned above that the aperture of the fracture is very thin, and
the fracture permeability is high.

\section{Discretizing fractures on the fine grid}

\subsection{Discrete Fracture Model (DFM)}
\label{sec:DFN}

Our first approach is based on representing
the fractures on the fine grid as the edges of finite element mesh.
This allows meshing the fractures more accurately; however, it can
be expensive when the fracture distribution is irregular.
Following a standard convention, we call ``matrix'' the region
that excludes all the fractures.
In the DFM (see e.g.,  \cite{Baca84}) as implemented in this paper,
we assume that the permeability does not vary along
the fractures, which coincide with the edges (or faces) of
finite element mesh. The fracture aperture
is taken into account in the finite element discretization as
a lower dimensional lines or surfaces. More precisely,
the discretization of the fractures on low dimensional surfaces
are added to finite element discretization of the matrix system.
Below, we demonstrate this main idea in our two dimensional example.



\begin{figure}[htb]
  \centering
  \includegraphics[width=0.65 \textwidth]{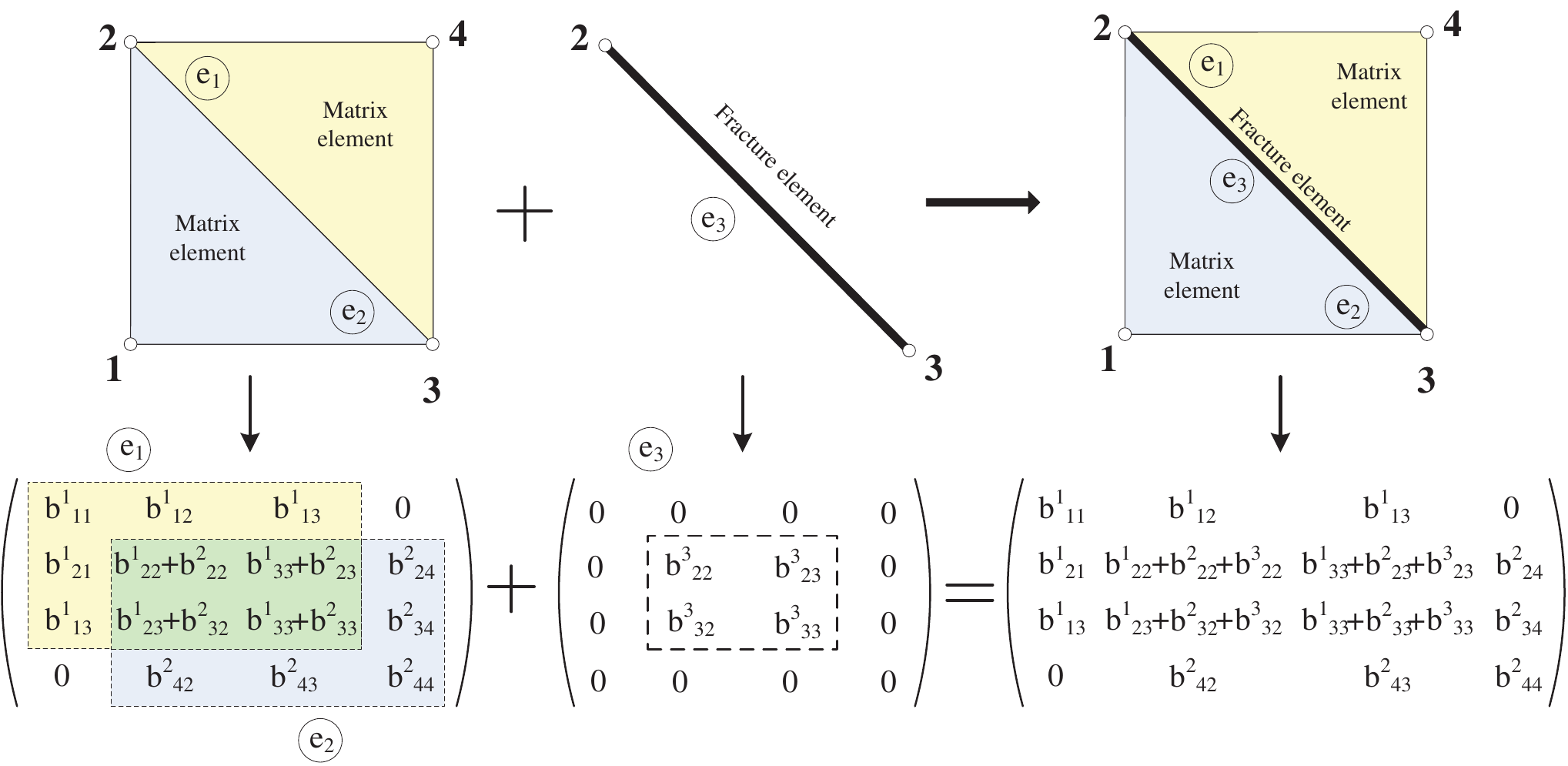}
  \caption{Illustration of DFM.}
  \label{fig:frame_int}
\end{figure}

Based on the assumptions above on the DFM, the fractures can be treated as one dimensional.
One-dimensional element is introduced in addition to the two-dimensional element for the discretization of the matrix.
Thus the Eqn. \eqref{eq:original} will be discretized in two-dimensional form for the matrix and in one-dimensional form for the fractures. The whole domain $D$ can be represented by
\begin{equation}
D=D_m\cup(\cup_i\, D_{\text{frac},i}),
\end{equation}
where the subscript $m$ and $\text{frac}$ represent the matrix and the fracture regions, respectively and $i$ refers to $i$th fracture.
Note that $D_m$ is a two-dimensional domain and $D_{\text{frac},i}$ is a one-dimensional domain.
We write the finite element discretization corresponding to
 Eqn. \eqref{eqn:fine-scale prb} as (see also Figure \ref{fig:frame_int}
for illustration)
%
\begin{equation} \label{eqn:integral}
\int_D \kappa(x) \nabla u \cdot \nabla v \, dx=
\int_{D_m} \kappa(x) \nabla u \cdot \nabla v \, dx
+\sum\limits_i \int_{D_{\text{frac},i}} \kappa(x) \nabla u \cdot \nabla v \, dx= \int_D f v \, dx.
\end{equation}
%

\subsection{Embedded Fracture Model (EFM)}

In this section, we present a fine-discretization of fracture network
following EFM proposed by Lee et al. (\cite{Lee01}).
EFM allows avoiding complex meshing for fracture network.
We would like again to emphasize that our goal is not to compare DFM with EFM, but rather
to develop GMsFEM framework which uses both fine-grid discretization techniques.
In EFM, long fractures are treated as low dimensional
objects and their interaction with the matrix is modeled separately.
In our two dimensional simulations, we
treat the long fractures separately as one-dimensional problems.
Eqn. \eqref{eq:original} can be written as
\begin{align}
\label{eq:hierarchical_model_1}
-\mbox{div} \big( \kappa(x) \, \nabla u  \big) +\Psi ^{mf}=f \quad \text{in} \quad D,\\
\label{eq:hierarchical_model_2}
-\mbox{div} \big( \kappa(x) \, \nabla u  \big) +\Psi ^{fm}=f \quad \text{in} \quad \text{fractures},
\end{align}
where $\Psi ^{mf}$ and $\Psi ^{fm}$ represent the
interaction between the long fractures and the rest of the system. Those two
terms can be calculated following \cite{Lee01}.
We denote
the linear system corresponding to the discrete form of the problem above is
\begin{align}\label{eq:hierarchical_matrix}
\begin{bmatrix}
A_{m} & B_{mf}^{1}&\cdots&B_{mf}^{Nf}\\
B_{fm}^{1}&B^{1}&\ldots&0\\
\vdots & \vdots & \ddots & \vdots\\
B_{fm}^{Nf}&0&\ldots&B^{Nf}
\end{bmatrix}
\begin{bmatrix}
u^{m}\\u^{1}\\\vdots\\u^{Nf}
\end{bmatrix}
=\begin{bmatrix}
f^{m}\\f^{1}\\ \vdots\\f^{Nf}
\end{bmatrix}
.
\end{align}
Here, $Nf$ denotes the total number of long fractures, $A_{m}$ is the matrix obtained over the domain excluding the
long fractures. {$B_{mf}^{i}$ denotes the interaction of the
$i${th} long fracture with the matrix and $ B_{fm}^{i}$ denotes the matrix interaction with the
$i$th long fracture, where $i=1,2, \cdots, Nf$.}
Here, we use EFM to handle long fractures, though
it can be used for shorter fractures.

\section{GMsFEM}\label{sec:GMsFEM}
In this section, we will give a brief description of the GMsFEM
for heterogeneous flow problems.
More details can be found in \cite{egh12, eglp13}. In the following,
we give a general outline of the GMsFEM. Then, we will discuss
the multiscale basis construction.

\noindent
{\bf Offline computations:}
\begin{itemize}
\item[Step 1] Coarse grid generation.
\item[Step 2] Construction of snapshot space that will be used to compute an offline space.
\item[Step 3] Construction of a small dimensional offline space by performing dimension reduction in the space of local snapshots.
\end{itemize}

If the problem has a parameter, one can use the offline space to construct
multiscale basis functions in the online stage.
In the above outline, the offline space can be reused if we
change the input of the model.
Given the computational domain, a coarse grid can be constructed
and local problems are solved on coarse neighborhoods
to obtain the snapshot spaces.
Then, smaller dimensional offline spaces are obtained
from the snapshot spaces by dimension reduction
via some spectral problems.
\subsection{Local basis functions}
\label{locbasis}

We now present the construction
of the basis functions
and the corresponding spectral problems
for obtaining a model reduction.


In the offline computation, we first construct a snapshot space $V_{\text{snap}}^{\omega_i}$ for each coarse neighborhood $\omega_i$.
The snapshot space can be the space of all fine-scale basis functions
or the solutions of some local problems with various choices of boundary conditions.
In \cite{randomized2014},
 we use $\kappa$-harmonic extensions
to form a snapshot space.
For flows in fractured media, the snapshot vectors are computed in
the same way except that we use an appropriate fine-grid discretization taking into account
the fracture distribution. In specific,
given a fine-scale piecewise linear function $\delta_j^h(x)$ defined on
$\partial\omega_i$, we define $\psi_{j}^{\omega_i, \text{snap}}$ by
\begin{equation} \label{harmonic_ex}
-\text{div}(\kappa(x) \nabla \psi_{j}^{\omega_i, \text{snap}} ) = 0
\quad \text{in} \, \, \, \omega_i,
\end{equation}
where $\psi_{j}^{\omega_i, \text{snap}}=\delta_j^h(x)$ on $\partial\omega_i$.
{\it The local problem is solved taking into consideration the fracture distribution.}
Later, we use randomized boundary conditions to avoid computing
the full snapshot vectors and use only a few snapshot vectors.
In particular, we use DFM for computing all snapshot vectors. For EFM,
we follow a hierarchical approach and discretize  few fractures
using the interaction terms. We can also add them into the calculations
of snapshot vectors; however, in the paper, we consider EFM for handling
a few long fractures.

For brevity of notation we now omit the superscript $\omega_i$, yet it is assumed throughout this section that the offline space computations are localized to respective coarse subdomains.
Let $l_i$ be the number of functions in the snapshot space in the region $\omega_i$, and
$$
V_{\text{snap}} = \text{span}\{ \psi_{j}^{ \text{snap}}:~~~ 1\leq j \leq l_i \},
$$
for each coarse subdomain $\omega_i$.
Denote
$$
R_{\text{snap}} = \left[ \psi_{1}^{\text{snap}}, \ldots, \psi_{l_i}^{\text{snap}} \right].
$$

In order to construct the offline space $V_{\text{off}}^\omega$, we perform a dimension reduction in the local snapshot space using an auxiliary spectral decomposition. The analysis in \cite{egw10} motivates the following
eigenvalue problem in the space of snapshots:
\begin{equation} \label{offeig}
A^{\text{off}} \Psi_k^{\text{off}} = \lambda_k^{\text{off}} S^{\text{off}} \Psi_k^{\text{off}},
\end{equation}
where
\begin{equation*}
 \displaystyle A^{\text{off}} = [a_{mn}^{\text{off}}] = \int_{\omega} \kappa(x) \nabla \psi_m^{\text{snap}} \cdot \nabla \psi_n^{\text{snap}} = R_{\text{snap}}^T A R_{\text{snap}}
 \end{equation*}
 \begin{center}
 and
 \end{center}
 \begin{equation*}
 \displaystyle S^{\text{off}} = [s_{mn}^{\text{off}}] = \int_\omega  \widetilde{\kappa}(x)\psi_m^{\text{snap}} \psi_n^{\text{snap}} = R_{\text{snap}}^T S R_{\text{snap}}.
\end{equation*}
The above integrals take into account the fracture distributions
(see Eqn. \eqref{eqn:integral}). We present the details of $ \widetilde{\kappa}$ later.

Here $A$ and $S$ denote analogous fine-scale matrices as defined by
\begin{equation*}
A_{ij} = \int_{D} \kappa(x) \nabla \phi_i \cdot \nabla \phi_j \, dx
\quad
S_{ij} = \int_{D}  \widetilde{\kappa}(x)   \phi_i  \phi_j \, dx,
\end{equation*}
where $\phi_i$ is the fine-scale basis function, and $\widetilde{\kappa}(x)$ is defined in the next subsection.

To generate the offline space we then choose the smallest $M^{\omega}_{\text{off}}$ eigenvalues from Eqn.~\eqref{offeig} and form the corresponding eigenvectors in the space of snapshots by setting
$\psi_k^{\text{off}} = \sum_{j=1}^{l_i} \Psi_{kj}^{\text{off}} \psi_j^{\text{snap}}$ (for $k=1,\ldots, M^{\omega}_{\text{off}}$), where $\Psi_{kj}^{\text{off}}$ are the coordinates of the vector $\psi_{k}^{\text{off}}$.

\subsection{Global coupling}
\label{globcoupling}

In this subsection, we discuss the offline space and the variational formulation for a continuous Galerkin approximation of Eqn. \eqref{eq:original}. We begin with an initial coarse space $V^{\text{init}}_0 = \text{span}\{ \chi_i \}_{i=1}^{N}$.
{Denote} $N$ to be the number of coarse neighborhoods. Here, $\chi_i$ are the standard multiscale partition of unity functions defined by
\begin{eqnarray} \label{pou}
-\text{div} \left( \kappa(x) \, \nabla \chi_i  \right) = 0 \quad K \in \omega_i \\
\chi_i = g_i \quad \text{on} \, \, \, \partial K, \nonumber
\end{eqnarray}
for all $K \in \omega_i$, where $g_i$ is a continuous function on $\partial K$ and is linear on each edge of $\partial K$.
We note that pointwise energy $\widetilde{\kappa}$ required for the eigenvalue problems is defined as
\begin{equation*}
\widetilde{\kappa} = \kappa \sum_{i=1}^{N} H^2 | \nabla \chi_i |^2,
\end{equation*}
where $H$ denotes the coarse mesh size.

We then multiply the partition of unity functions by the eigenfunctions in the offline space $V_{\text{off}}^{\omega_i}$ to construct the resulting basis functions
\begin{equation} \label{cgbasis}
\psi_{i,k} = \chi_i \psi_k^{\omega_i, \text{off}} \quad \text{for} \, \, \,
1 \leq i \leq N \, \, \,  \text{and} \, \, \, 1 \leq k \leq M_{\text{off}}^{\omega_i},
\end{equation}
where $M_{\text{off}}^{\omega_i}$ denotes the number of offline eigenvectors that are chosen for each coarse node $i$.
We note that the construction in Eqn.~\eqref{cgbasis} yields continuous basis functions due to the multiplication of
 offline eigenvectors with the initial (continuous) partition of unity function. Next, we define the continuous Galerkin
 spectral multiscale space as
\begin{equation} \label{cgspace}
V_{\text{off}}  = \text{span} \{ \psi_{i,k} : \,  \, 1 \leq i \leq N \, \, \,  \text{and} \, \, \, 1 \leq k \leq M_{\text{off}}^{\omega_i}  \}.
\end{equation}
Using a single index notation, we may write $V_{\text{off}} = \text{span} \{ \psi_{i} \}_{i=1}^{N_c}$, where $N_c =\sum_{i=1}^{N}M_{\text{off}}^{\omega_{i}}$
denotes the total number of basis functions in the space $V_{\text{off}}$. We also construct an operator matrix $R_0^T = \left[ \psi_1 , \ldots, \psi_{N_c} \right]$ (where $\psi_i$ are used to denote the nodal values of each basis function defined on the fine grid), for later use in this subsection.

Below, we will display in detail the multiscale formulation  using DFM and EFM, respectively.
\subsubsection{Multiscale DFM approach}
We seek $u_{\text{ms}}(x) = \sum_i c_i \psi_i(x) \in V_{\text{off}}$ such that
\begin{equation} \label{cgvarform}
a(u_{\text{ms}}, v) = (f, v) \quad \text{for all} \,\,\, v \in V_{\text{off}}.
\end{equation}
Note that the offline space $V_{\text{off}}$ is an approximation of all the nodal basis, including the ones on the long fractures.
Therefore, the variational form in \eqref{cgvarform} yields the following linear algebraic system
\begin{equation}
A_0 U_0 = F_0,
\end{equation}
where $U_0$ denotes the nodal values of the discrete solution, and
\begin{equation*}
A_0 = [a_{IJ}] = \int_D \kappa(x)  \, \nabla \psi_I \cdot \nabla \psi_J \, dx \quad \text{and} \quad F_0 = [f_I] = \int_D f \psi_I \, dx.
\end{equation*}
Using the operator matrix $R_0^T$, we may write $A_0 = R_0 A R_0^T$ and $F_0 = R_0 F$, where $A$ and $F$ are the standard, fine-scale stiffness matrix and forcing vector corresponding to the form in Eqn. {\eqref{eqn:integral}}. We also note that the operator matrix may be analogously used in order to project coarse scale solutions onto the fine grid.

\subsubsection{Multiscale EFM approach}

Similar to the multiscale DFM approach proposed above, we seek $u_{\text{ms}}(x) = \sum_i c_i \psi_i(x) \in V_{\text{off}}$ such that
\begin{equation*}
a(u_{\text{ms}}, v) = (f, v) \quad \text{for all} \,\,\, v \in V_{\text{off}}.
\end{equation*}
However, the nodal basis on the long fractures are excluded from the offline space $V_{\text{off}}$, i.e., the fine-scale nodal basis
are used on the long fractures.
It follows from Eqn. \eqref{eq:hierarchical_matrix},
\begin{align}\label{eq:hierarchical_ms}
\begin{bmatrix}
R_0A_{m}R_0^T &R_0 B_{mf}^{1}&\cdots&R_0B_{mf}^{Nf}\\
B_{fm}^{1}R_0^T&B^{1}&\ldots&0\\
\vdots & \vdots & \ddots & \vdots\\
B_{fm}^{Nf}R_0^T&0&\ldots&B^{Nf}
\end{bmatrix}
\begin{bmatrix}
u^{m}_{ms}\\u^{1}\\\vdots\\u^{Nf}
\end{bmatrix}
=\begin{bmatrix}
R_0f^{m}\\f^{1}\\ \vdots\\f^{Nf}
\end{bmatrix}
.
\end{align}
Note that the operator matrix may be analogously used in order to project coarse scale solutions onto the fine grid.



\section{Numerical result}
\label{sec:numerical}

In this section, we will present several numerical experiments to show
the performance of multiscale DFM and multiscale EFM approaches in the fracture modeling.
The simulation results using multiscale DFM is shown in Subsection \ref{subsec:dfm}. In Subsection \ref{subsec:efm},
a few tests are conducted to verify the performance of multiscale EFM approach. Moreover, we combine those two approaches
in the modeling of short and long fractures and list the results in Subsection \ref{subsec:dfm_efm}. Further, the relation between the basis selection and the types of fractures is studied in Subsection \ref{subsec:adaptive}.
We also present a randomized snapshot calculations to reduce the
computational cost associated with computing the snapshot space.

%
We take the domain {$D$} to be a square,
set the forcing term $f=0$ and impose a bilinear boundary condition for the problem \eqref{eq:original}.
In our numerical simulations, we use a $10 \times 10$ coarse grid,
and each coarse grid block is further divided into $10\times 10$ fine grid blocks.
Thus, the whole computational domain is partitioned into a $100 \times 100$ fine grid.
The fine-scale solution is obtained
by discretizing \eqref{eq:original} following DFM with piecewise bilinear elements
on the fine grid and linear one-dimensional elements over the fractures.

 We recall that
$V_{\text{off}}$ denotes the offline space; $u$, $u_{\text{snap}}$ and $u_{\text{off}}$ denote the fine-scale, snapshot and offline solutions, respectively.
In the tables below, we will
compute the error $u-u_{\text{off}}$ using
the $L^2$ relative error and the energy relative error, which are defined as
\begin{equation}
\| u-u_{\text{off}}\|_{L^2_{\kappa}(D)} := \frac{ \| u - u_{\text{off}} \|_{L^2(V)} }{ \| u \|_{L^2(V)} }, \quad\quad
\| u-u_{\text{off}}\|_{H^1_{\kappa}(D)} := \frac{ a(u-u_{\text{off}}, u-u_{\text{off}})^{\frac{1}{2}} }{ a(u,u)^{\frac{1}{2}} },
\end{equation}
where the weighted $L^2$-norm is defined as $\|u\|_{L^2(V)} = \| {\kappa}^{\frac{1}{2}} u \|_{L^2(D)}$.
We will also compute the error $u_{\text{snap}}-u_{\text{off}}$ using the same norms
\begin{equation}
\| u_{\text{snap}}-u_{\text{off}}\|_{L^2_{\kappa}(D)} := \frac{ \| u_{\text{snap}} - u_{\text{off}} \|_{L^2(V)} }{ \| u_{\text{snap}} \|_{L^2(V)} }, \quad\quad
\| u_{\text{snap}}-u_{\text{off}}\|_{H^1_{\kappa}(D)} := \frac{ a(u_{\text{snap}}-u_{\text{off}}, u-u_{\text{off}})^{\frac{1}{2}} }{ a(u_{\text{snap}},u_{\text{snap}})^{\frac{1}{2}} }.
\end{equation}

\subsection{Numerical results with DFM}
\label{subsec:dfm}

\begin{figure}[htb]
\centering
 \subfigure[Fracture 1]{\label{fig:frac_various_1}
    \includegraphics[width = 0.30\textwidth, keepaspectratio = true]{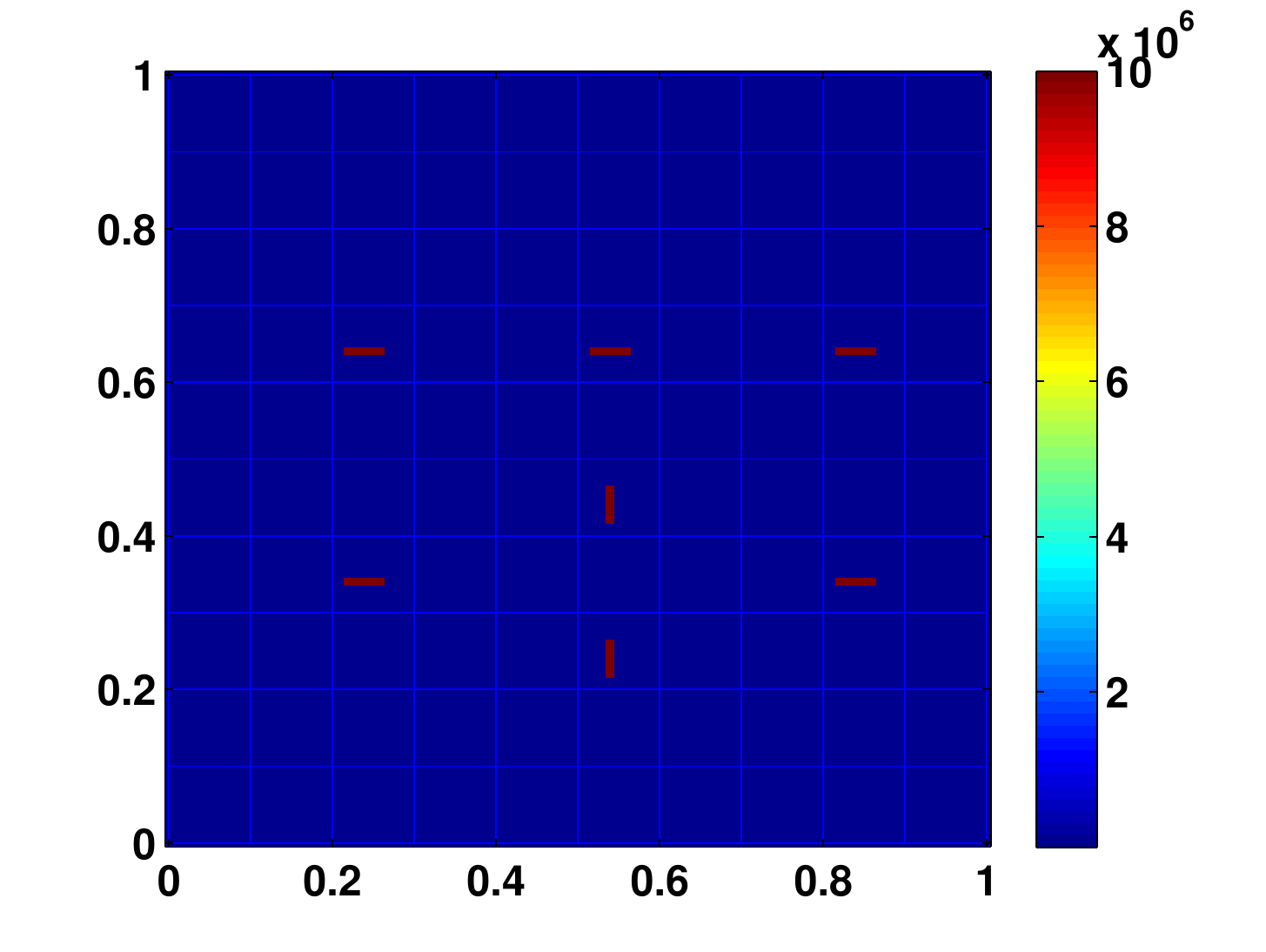}
   }
  \subfigure[Fracture 2]{\label{fig:frac_various_exp2}
    \includegraphics[width = 0.30\textwidth, keepaspectratio = true]{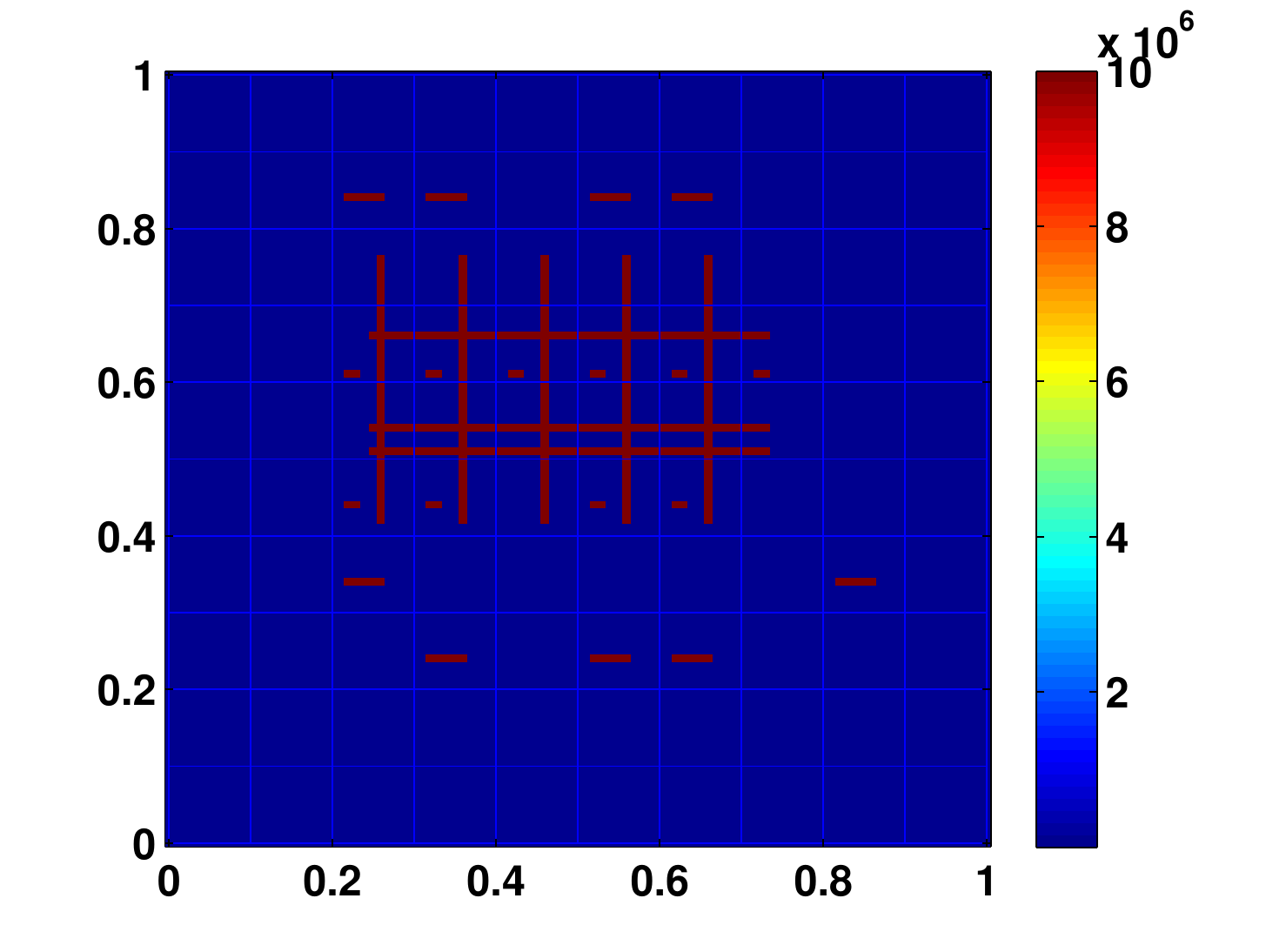}
   }
  \subfigure[Fracture 3]{\label{fig:frac_various_exp3}
    \includegraphics[width = 0.30\textwidth, keepaspectratio = true]{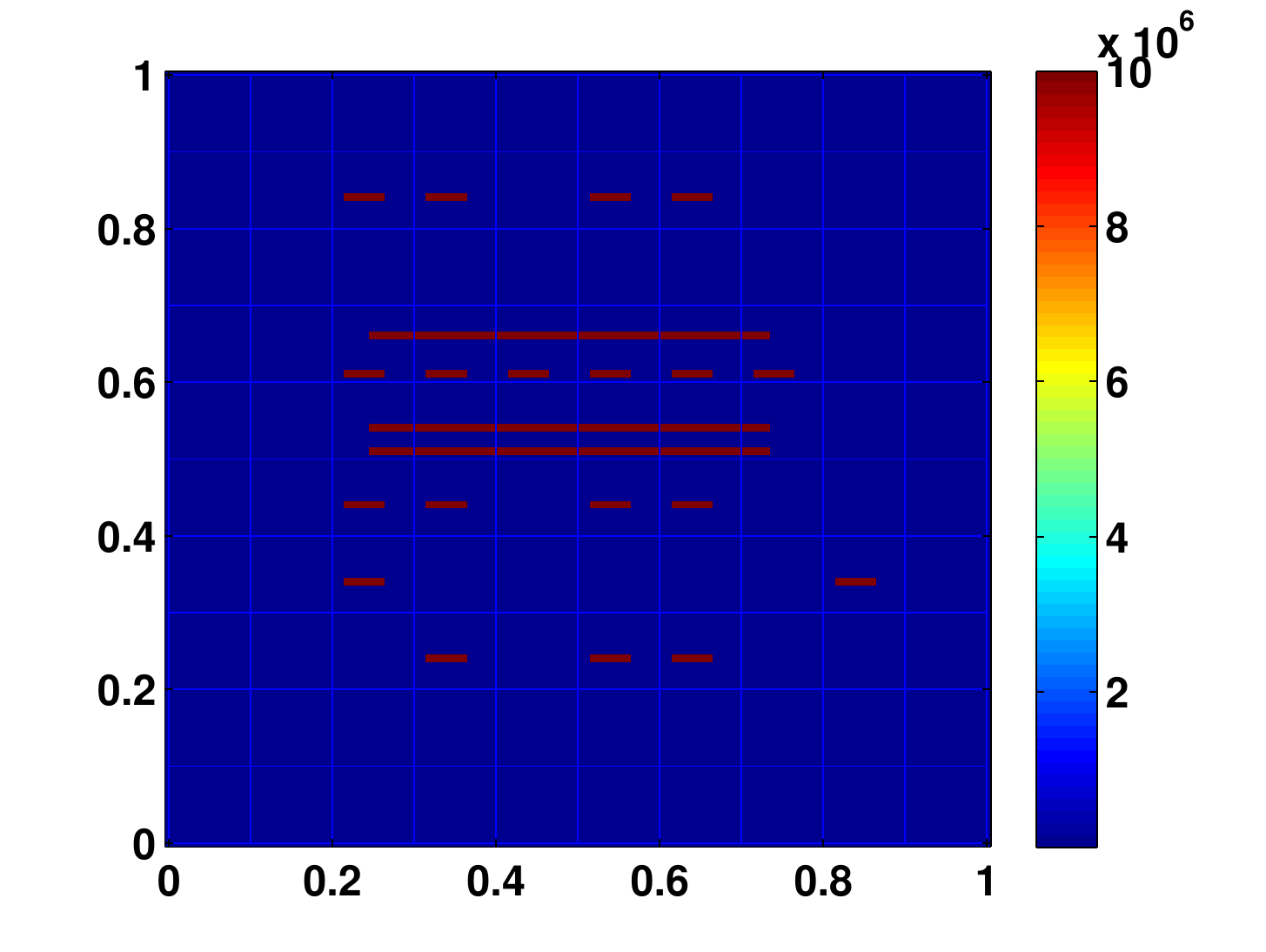}
   }

 \caption{Permeability fields with fractures. Fractures are shown with red, while the background is blue.}
 \label{fig:perms}
\end{figure}

\begin{figure}[htb]
\centering
 \subfigure[Fine-scale solution 1]{\label{fig:FEMsoln_frac_1}
    \includegraphics[width = 0.30\textwidth, keepaspectratio = true]{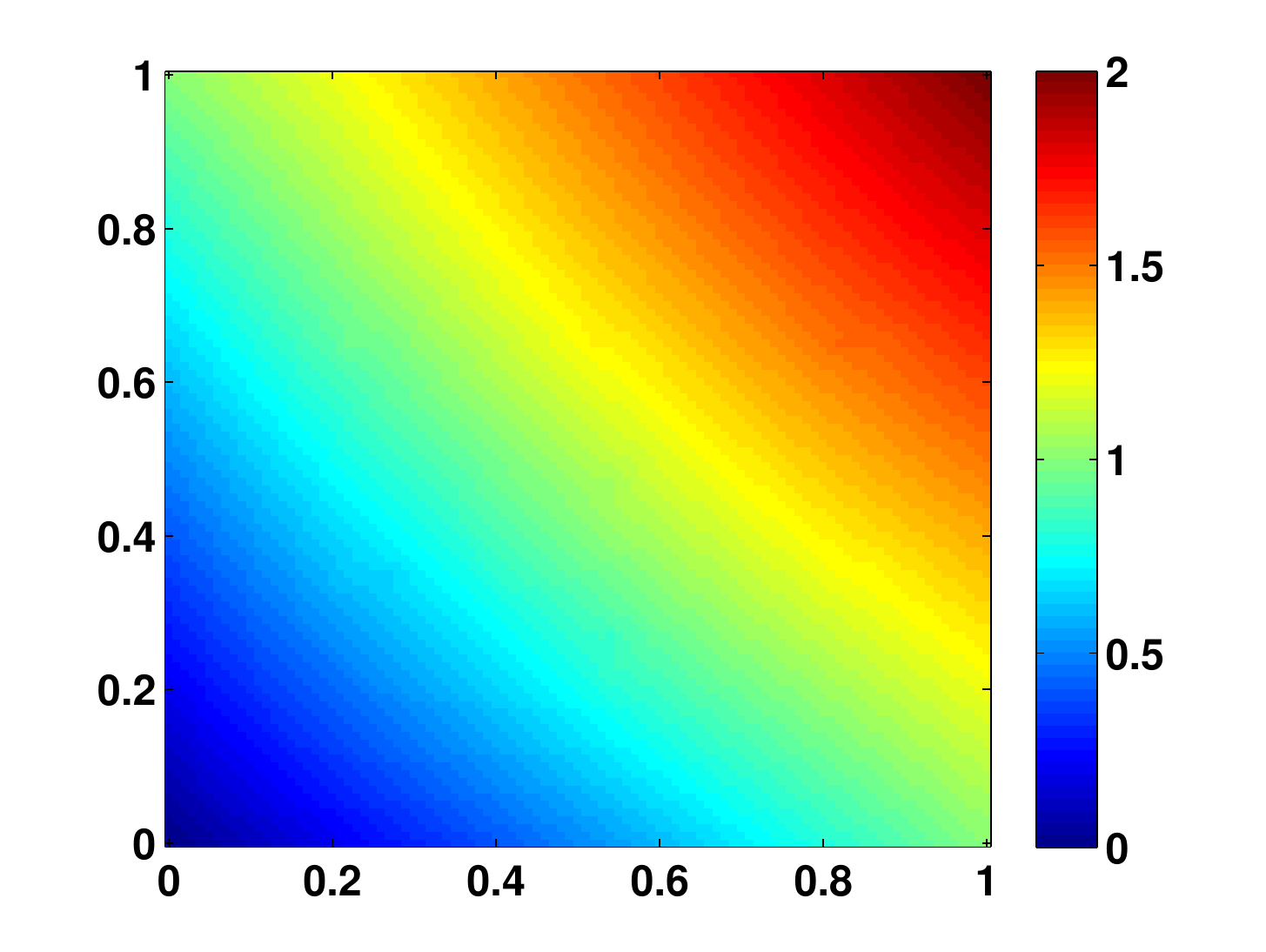}}
  \subfigure[Fine-scale solution 2]{\label{fig:FEMsoln_various_exp2}
    \includegraphics[width = 0.30\textwidth, keepaspectratio = true]{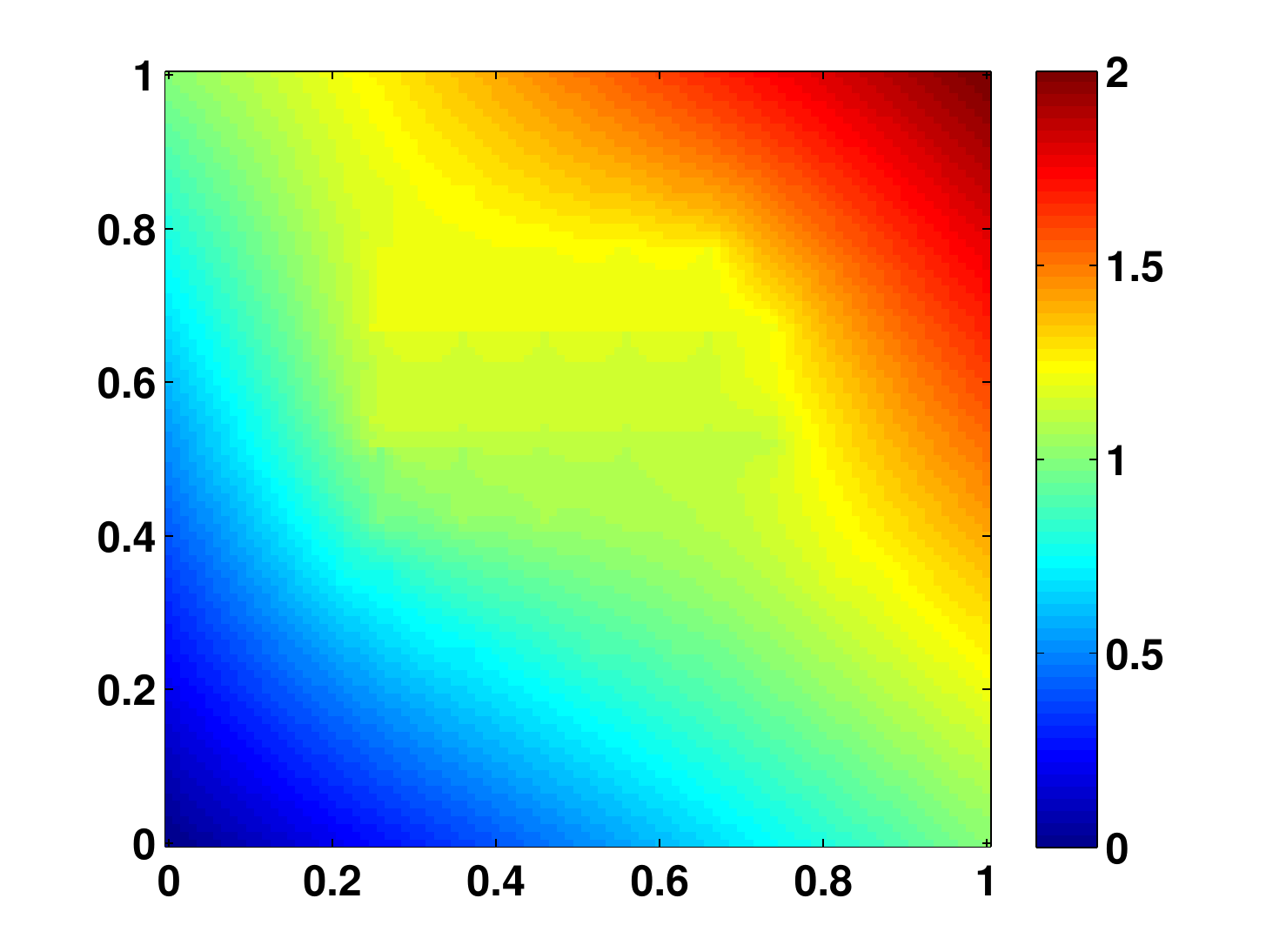}
   }
  \subfigure[Fine-scale solution 3]{\label{fig:FEMsoln_various_exp3}
    \includegraphics[width = 0.30\textwidth, keepaspectratio = true]{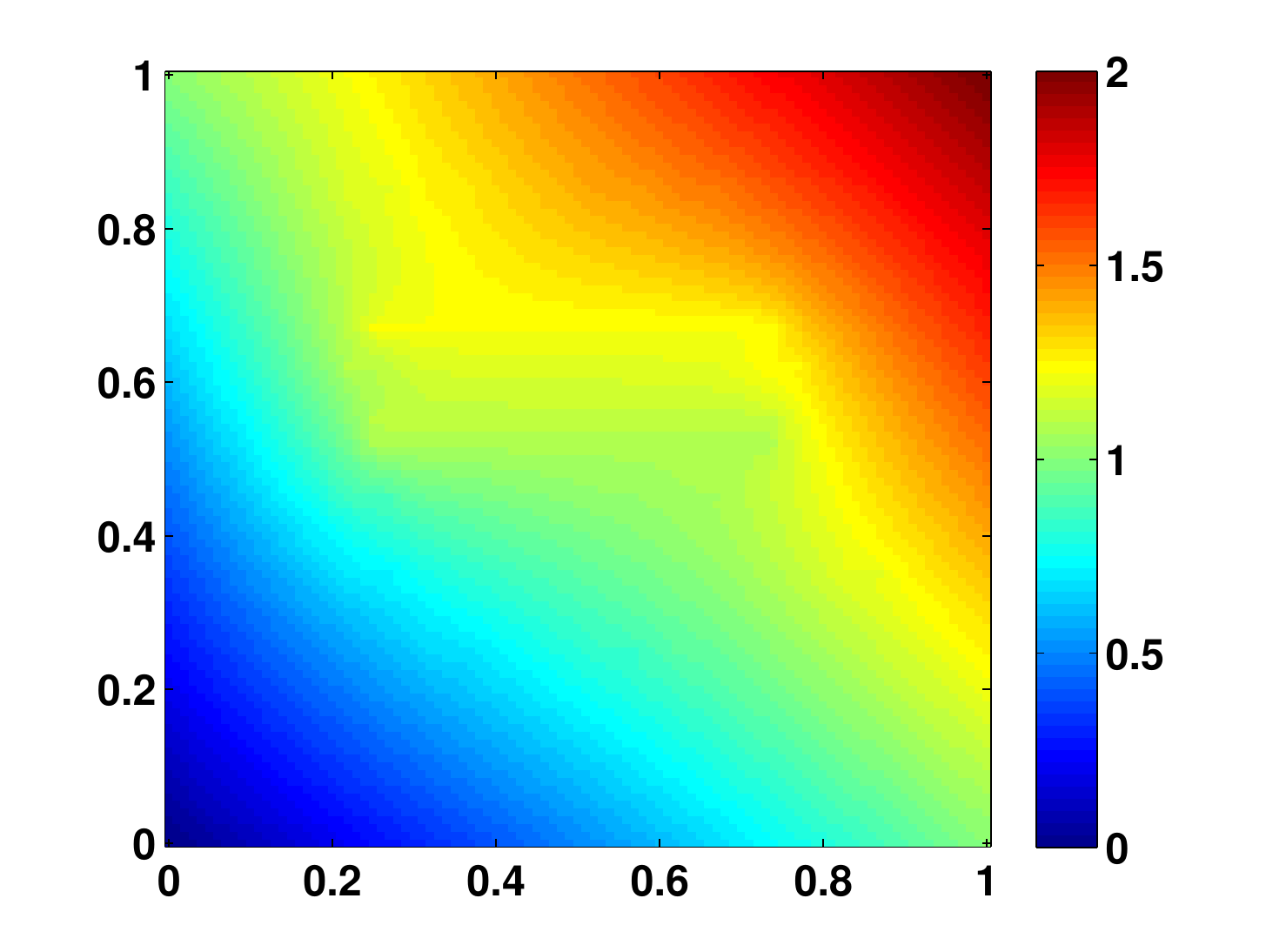}
   }
 \caption{Fine-scale solutions corresponding to the permeability fields in Fig. \ref{fig:perms}, respectively.}
 \label{fig:FEMsoln}
\end{figure}

\begin{figure}[htb]
\centering
 \subfigure[Coarse-scale solution 1]{\label{fig:umsfineSS_frac_1}
    \includegraphics[width = 0.30\textwidth, keepaspectratio = true]{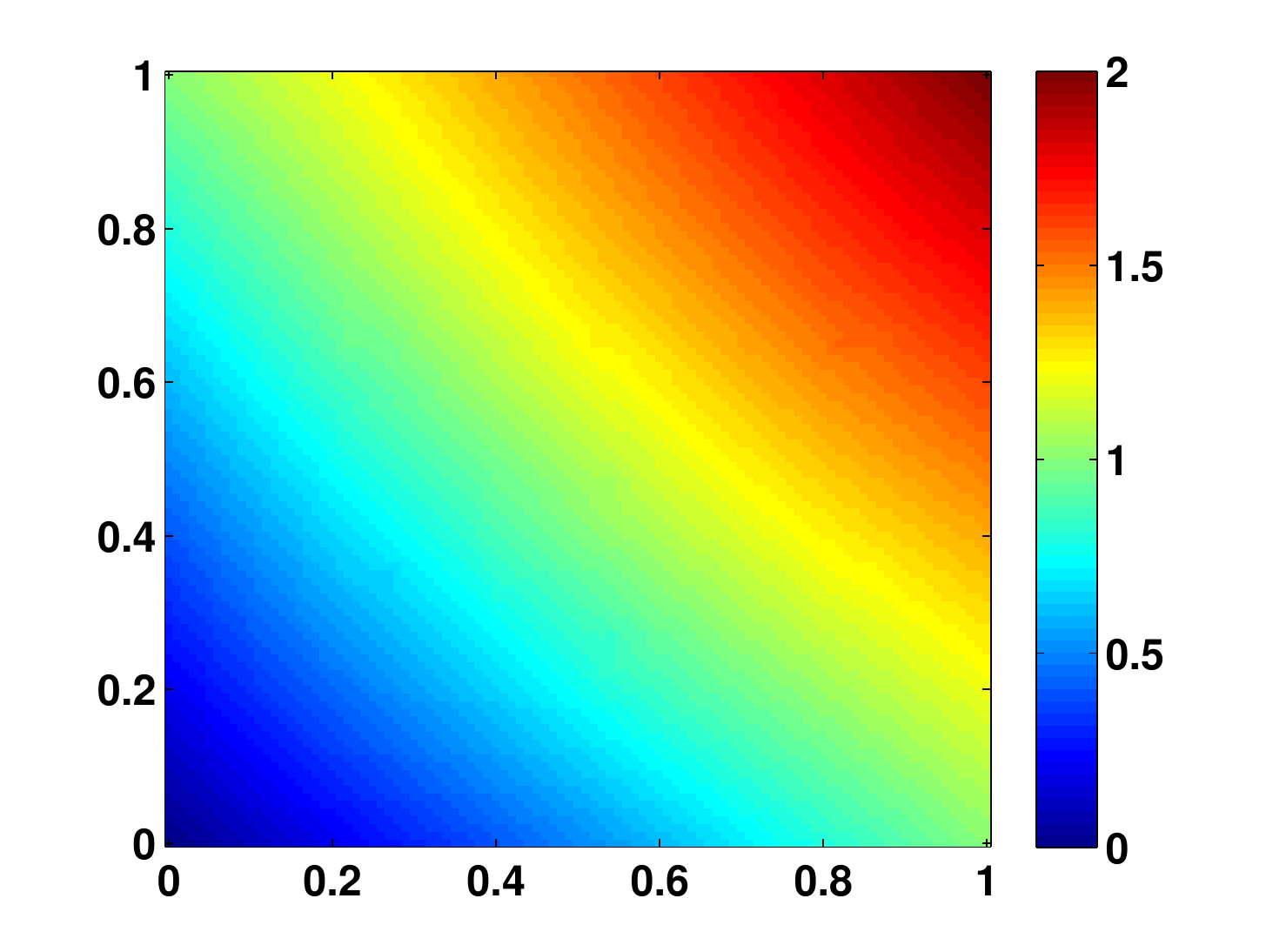}
   }
  \subfigure[Coarse-scale solution 2]{\label{fig:umsfineSS_various_exp2}
    \includegraphics[width = 0.30\textwidth, keepaspectratio = true]{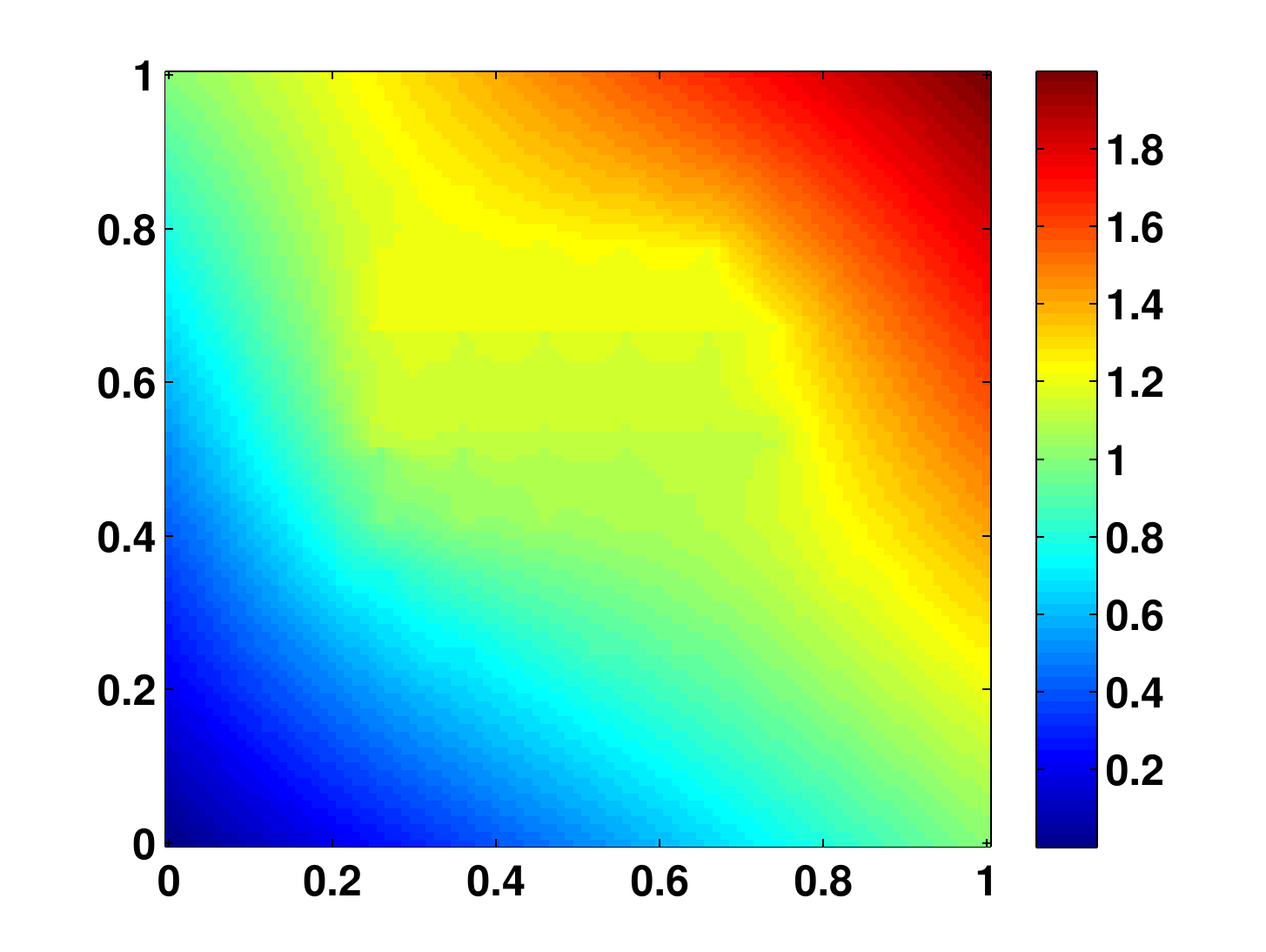}
   }
  \subfigure[Coarse-scale solution 3]{\label{fig:umsfineSS_various_exp3}
    \includegraphics[width = 0.30\textwidth, keepaspectratio = true]{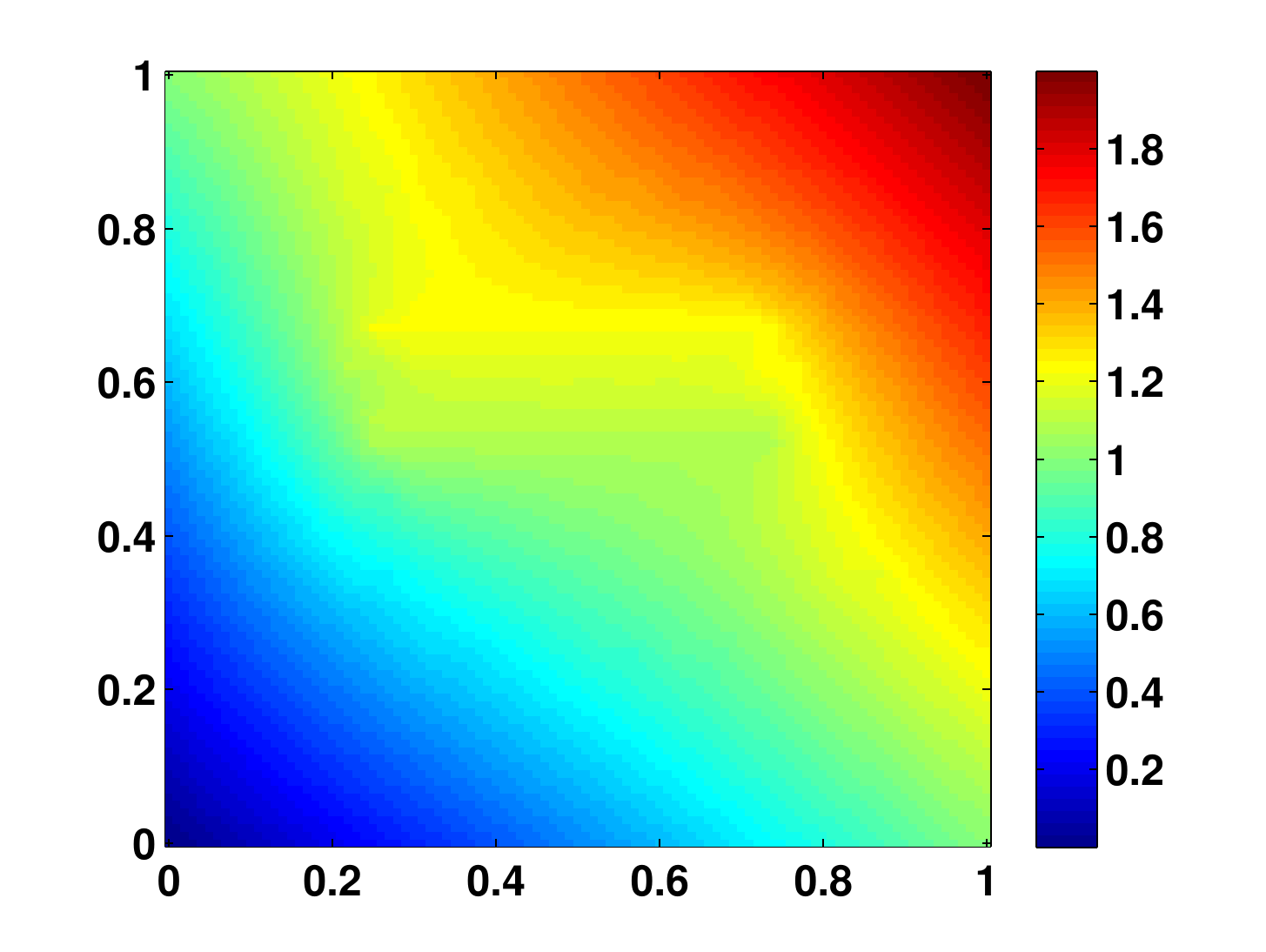}
   }
 \caption{Coarse-scale solutions corresponding to the permeability fields in Fig. \ref{fig:perms}, respectively}
 \label{fig:umsfineSS}
\end{figure}

In Fig. \ref{fig:perms}, we depict three fracture fields used in the simulations below.
In Fig. \ref{fig:frac_various_1}, there is at most one fracture in each coarse block without crossing the coarse edges, while in Figs. \ref{fig:frac_various_exp2} and \ref{fig:frac_various_exp3}, the fractures are
 more complicated and intersect the edges of coarse blocks.
The simulation result with fractures in Fig. \ref{fig:frac_various_1} with 1 basis per coarse node is very good (with the energy error of $2.52\%$) because
their effects can be localized. The simulation results with the other two fracture systems are shown in Tables \ref{table:frac_various_exp2} and \ref{table:frac_various_exp3}, respectively. We notice that one basis function per node
does not give a satisfactory for more complex fracture systems.
We also show the fine-scale solutions and the snapshot solutions (i.e. the coarse-scale solution corresponding to the largest offline space we select) corresponding to fractures in Fig. \ref{fig:perms}, which are shown in Figs. \ref{fig:FEMsoln} and \ref{fig:umsfineSS}, respectively.

\begin{figure}[htb]
\centering
 \subfigure[Fracture 1]{\label{fig:fracture_construction}
    \includegraphics[width = 0.30\textwidth, keepaspectratio = true]{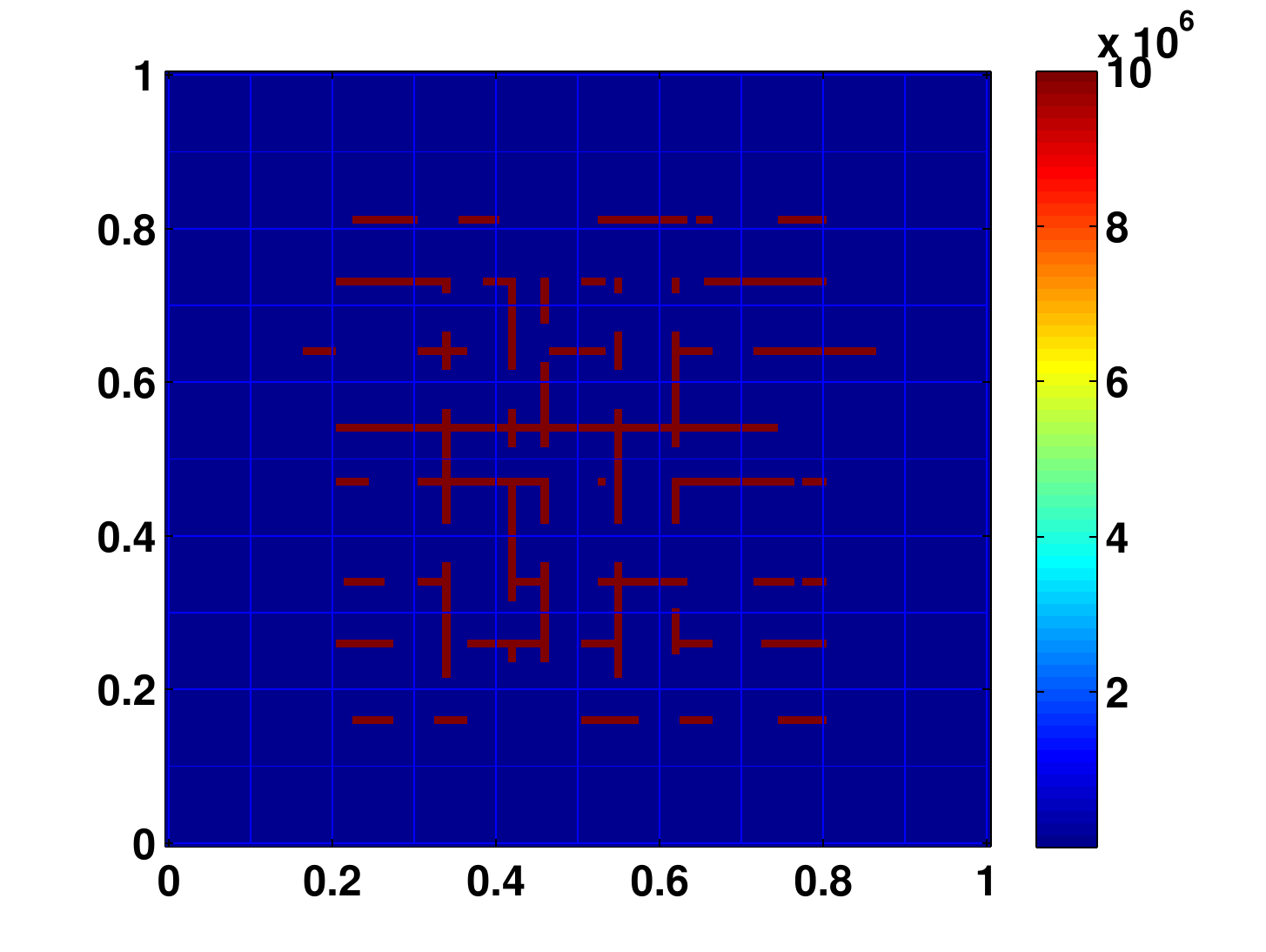}
   }
  \subfigure[Fracture 2]{\label{fig:fracture_construction1}
    \includegraphics[width = 0.30\textwidth, keepaspectratio = true]{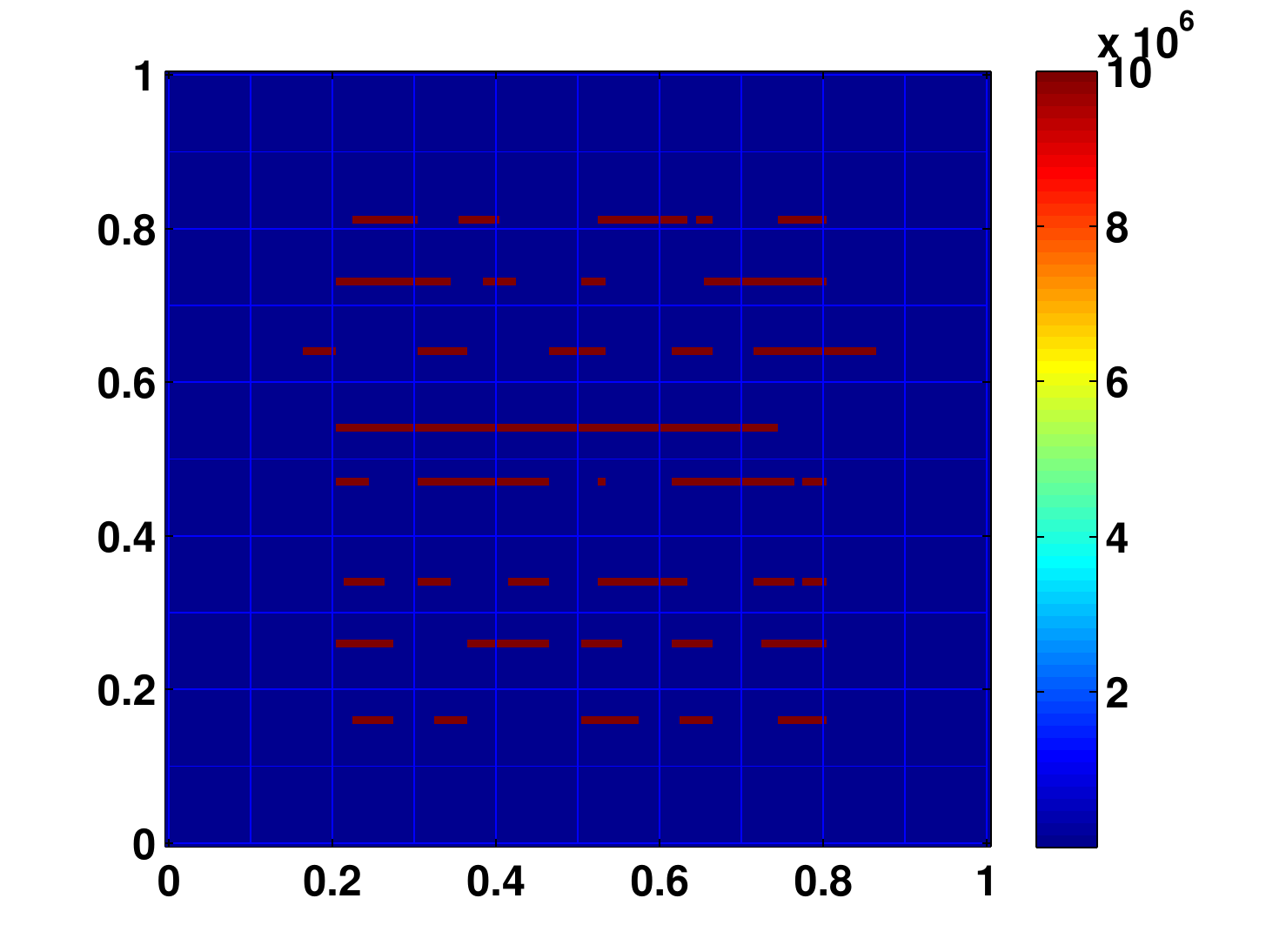}
   }
\subfigure[Fracture 3]{\label{fig:frac_tri}
    \includegraphics[width = 0.30\textwidth, keepaspectratio = true]{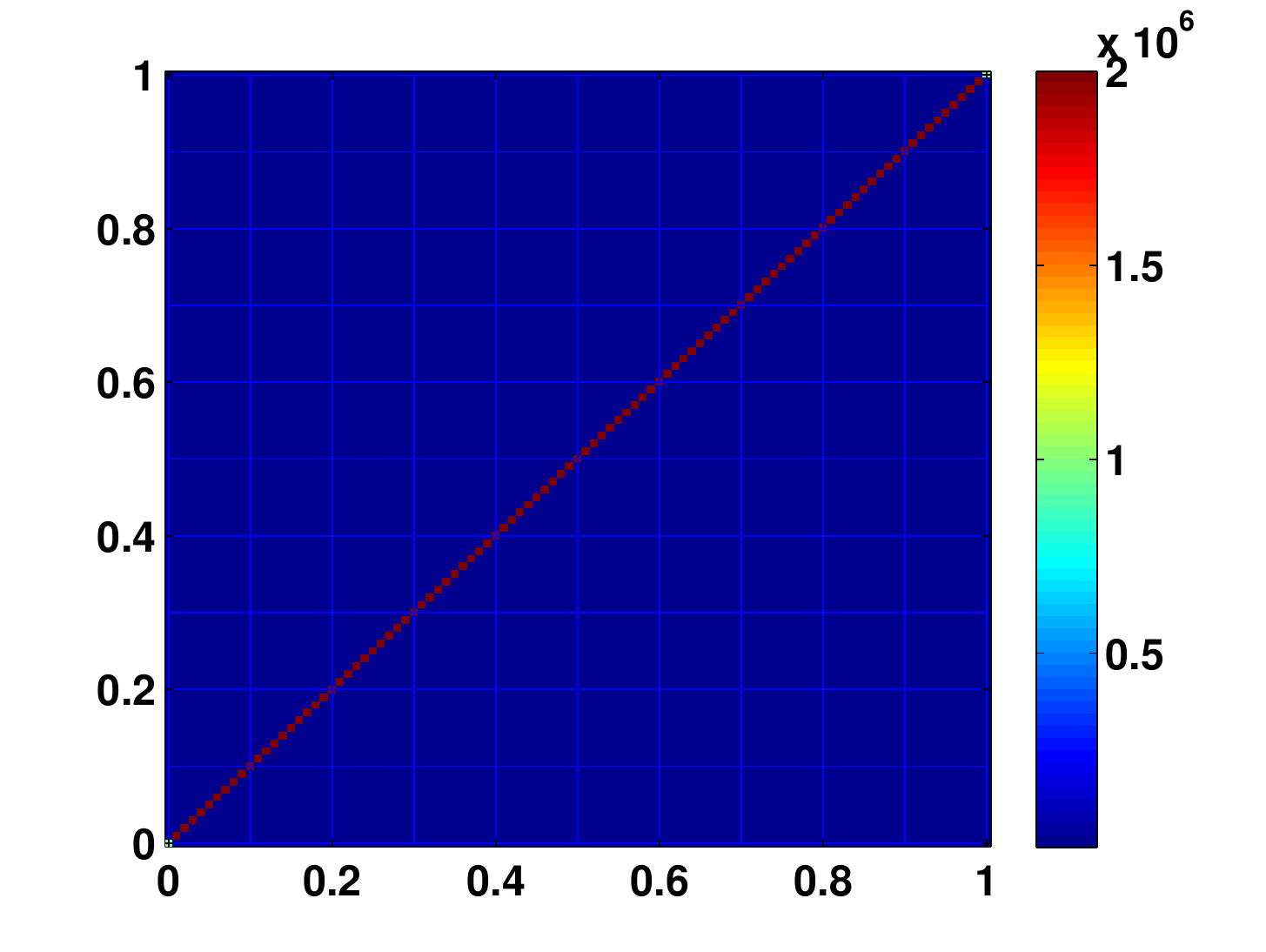}
   }
 \caption{Permeability fields with fractures. Fractures are shown with red, while the matrix is the blue region.}
 \label{fig:perms_add}
\end{figure}

\begin{figure}[htb]
\centering
  \subfigure[x-component of velocity]{\label{fig:x_direction_velocity}
    \includegraphics[width = 0.45\textwidth, height = 0.3\textwidth]{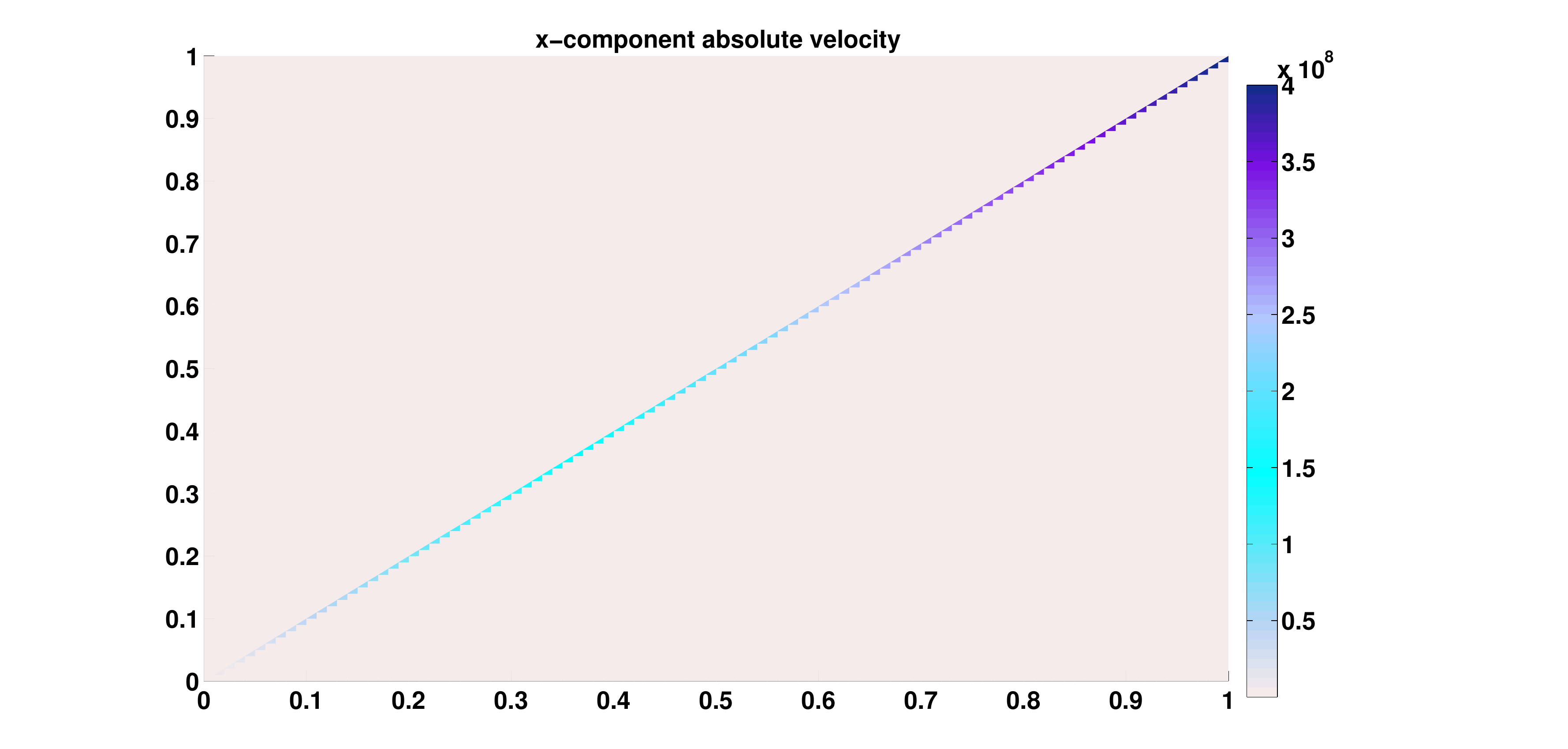}
   }
\subfigure[y-component of velocity]{\label{fig:y_direction_velocity}
    \includegraphics[width = 0.45\textwidth,height = 0.3\textwidth]{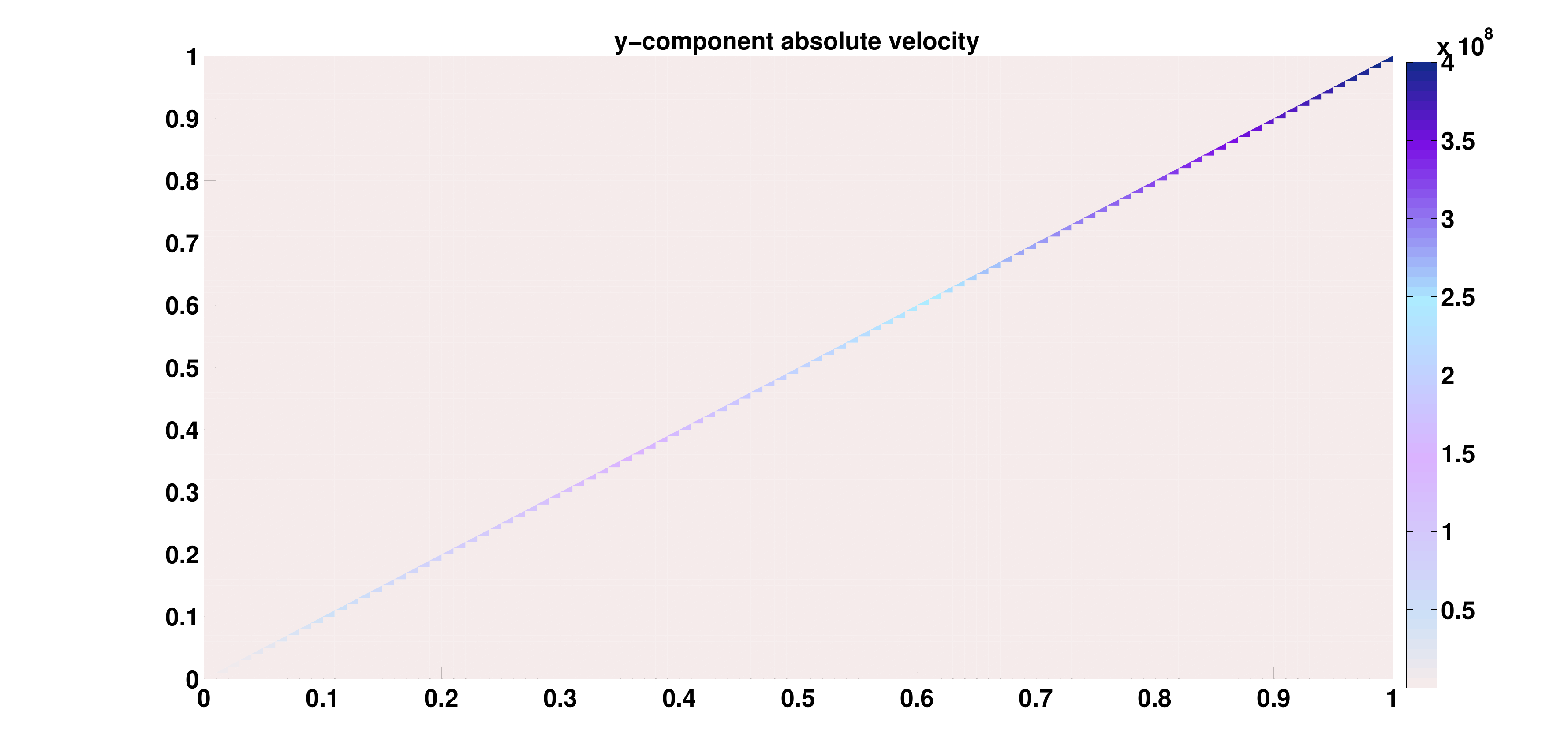}
   }
 \caption{Velocity field corresponding to the permeability field in Fig. \ref{fig:frac_tri}.}
 \label{fig:tri_velocity}
\end{figure}

The convergence history of fracture fields in Figs. \ref{fig:frac_various_exp2} and \ref{fig:frac_various_exp3} are displayed in Tables \ref{table:frac_various_exp2} and \ref{table:frac_various_exp3}. We observe that the offline solution will converge to the fine-scale solution as we enrich
the offline space in each coarse neighborhood. The energy error decreases from $27.25\%$ to $7.4\%$ as we add 4 more local
 basis in each coarse neighborhood from Table \ref{table:frac_various_exp2}. In Tables \ref{table:fracture_construction}
 and \ref{table:fracture_construction1}, we show the numerical results for fracture fields in Figs.
 \ref{fig:fracture_construction} and \ref{fig:fracture_construction1}.
%
\begin{figure}[htb]
\centering
 \subfigure[Fine-scale solution 1]{\label{fig:FEMsoln_fracture_construction}
    \includegraphics[width = 0.30\textwidth, keepaspectratio = true]{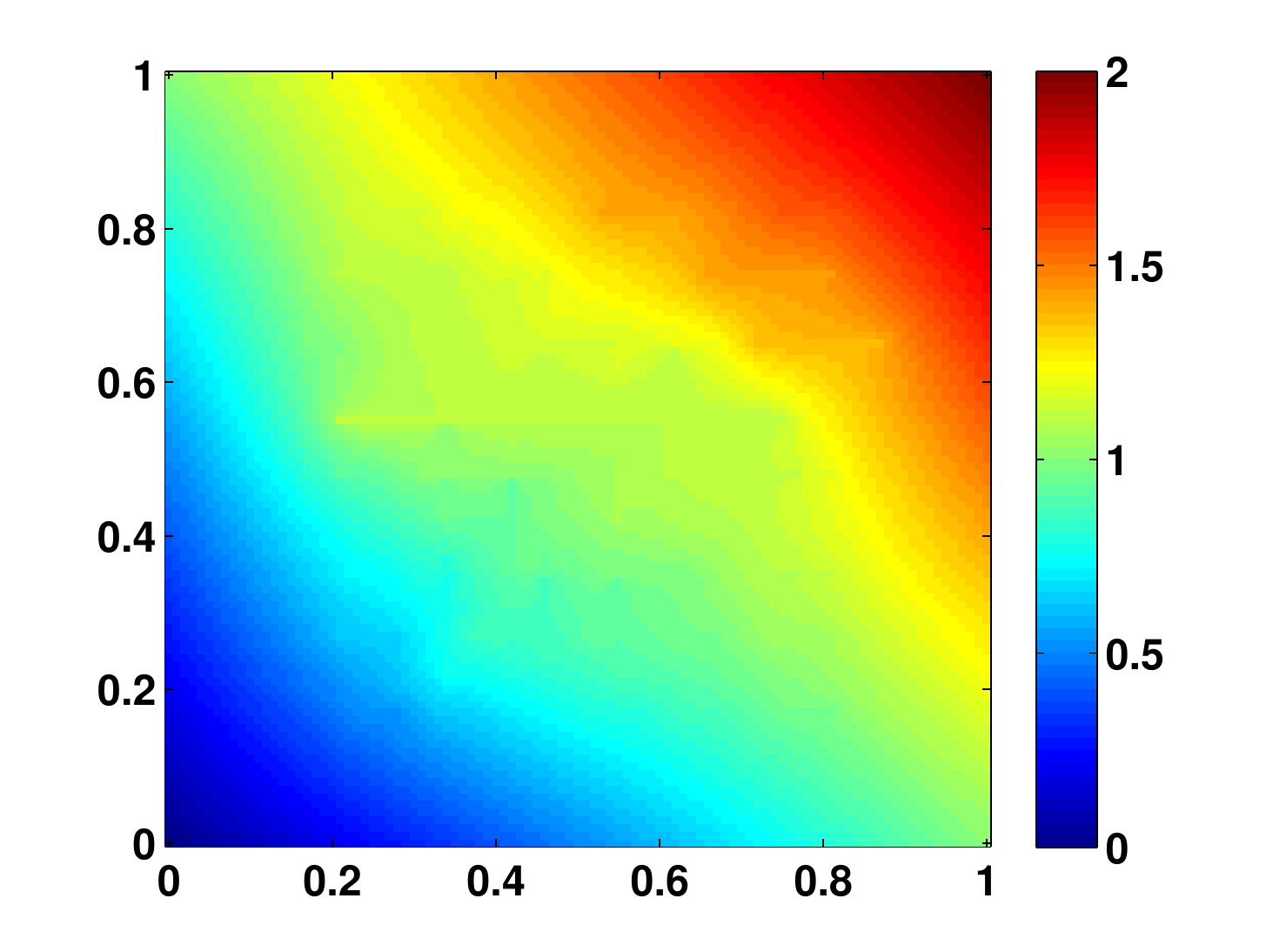}
   }
  \subfigure[Fine-scale solution 2]{\label{fig:FEMsoln_fracture_construction1}
    \includegraphics[width = 0.30\textwidth, keepaspectratio = true]{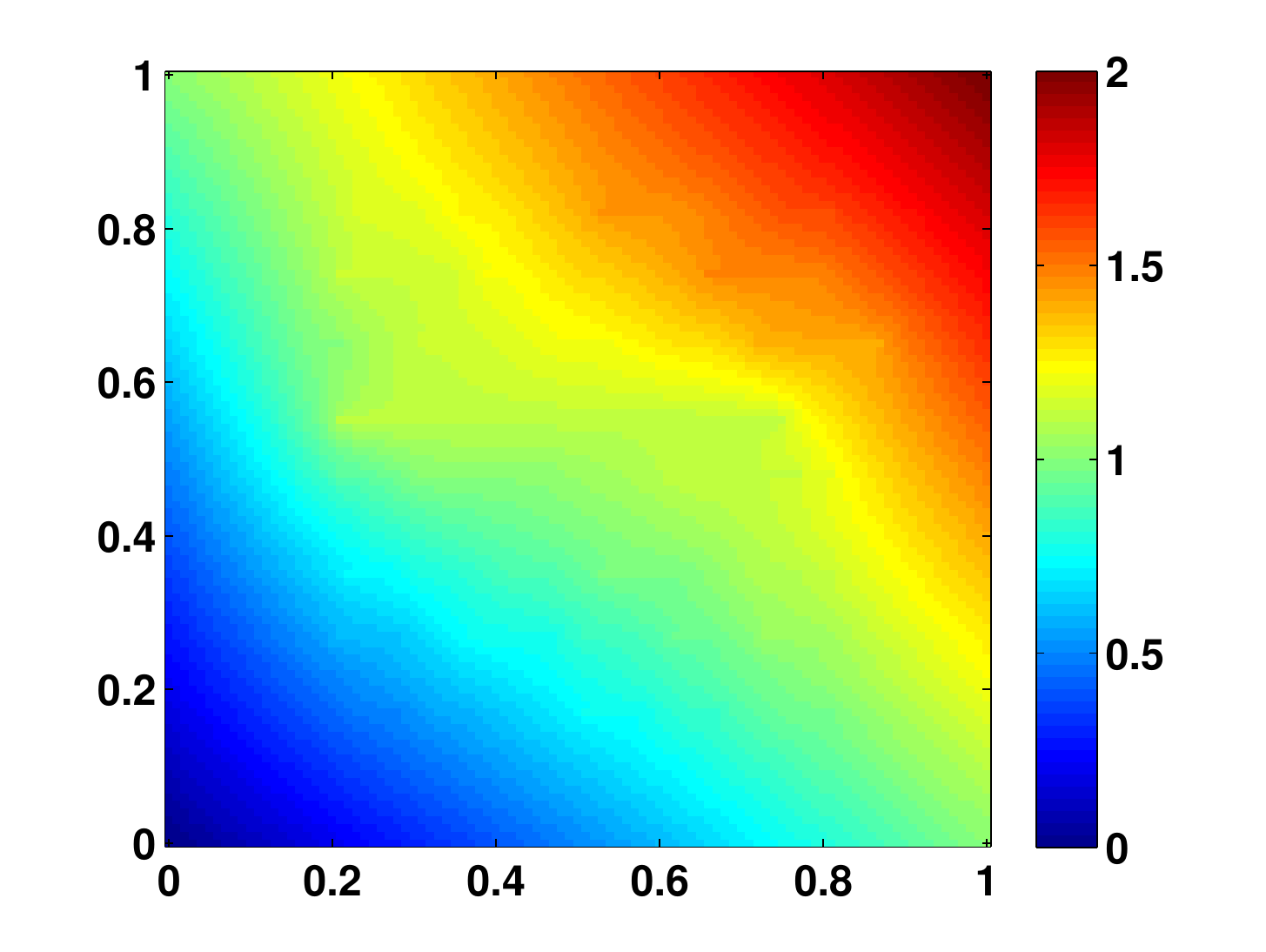}
   }
 \subfigure[pressure solution]{\label{fig:diag_solution}
    \includegraphics[width = 0.30\textwidth, keepaspectratio = true]{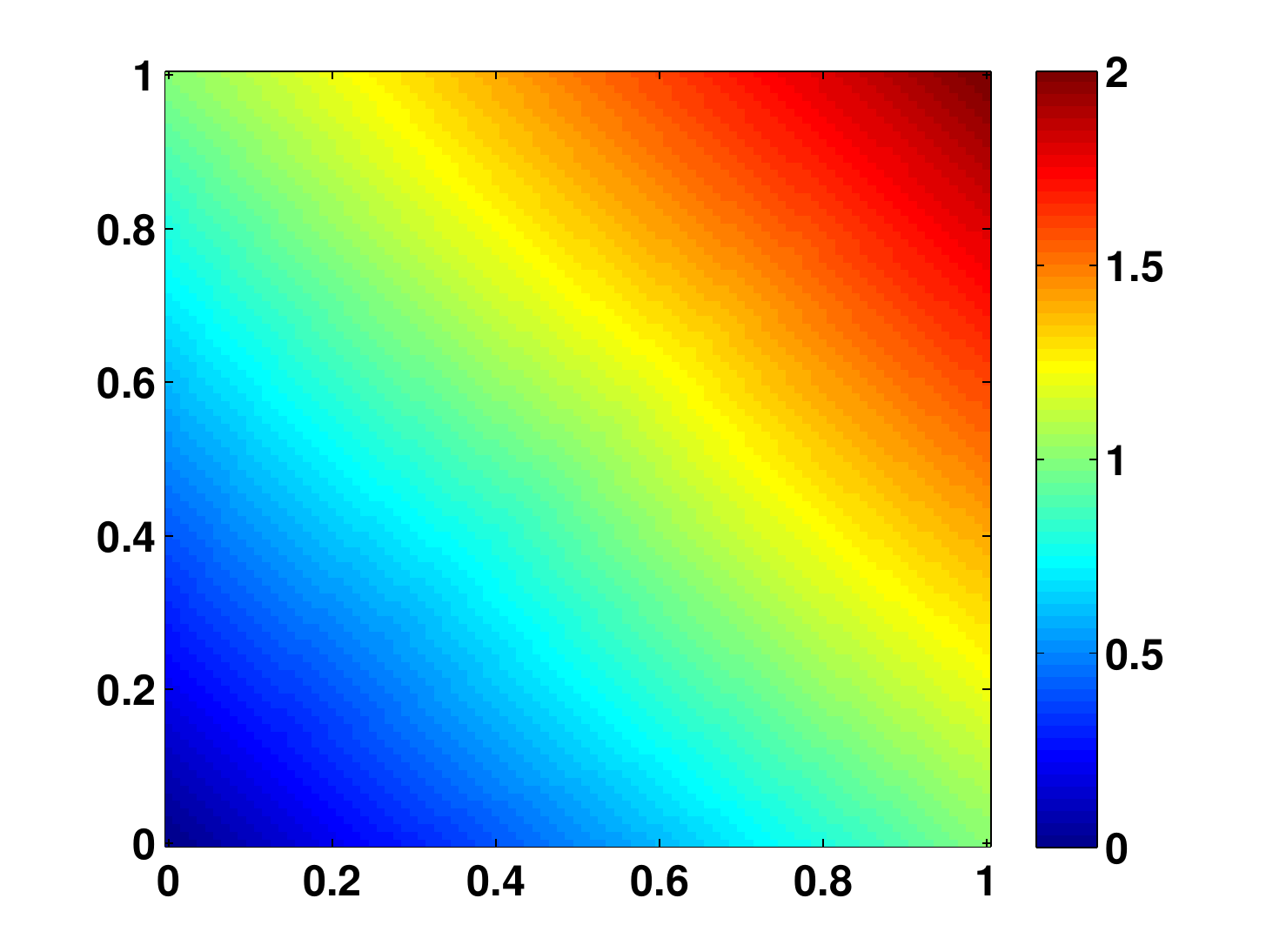}
   }
 \caption{Fine-scale solutions corresponding to the permeability fields in Fig. \ref{fig:perms_add}, respectively.}
 \label{fig:FEMsoln_add}
\end{figure}
\begin{table}[htb]
\centering
\begin{tabular}{|c|c|c|c|c|c|c|}
\hline
\multirow{2}{*}{$\text{dim}(V_{\text{off}})$} &
\multicolumn{2}{c|}{  $\|u-u_{\text{off}} \|$ (\%) } &
\multicolumn{2}{c|}{  $\|u_{\text{snap}}-u_{\text{off}} \|$ (\%) }\\
\cline{2-5} {}&
$\hspace*{0.8cm}   L^{2}_\kappa(D)   \hspace*{0.8cm}$ &
$\hspace*{0.8cm}   H^{1}_\kappa(D)  \hspace*{0.8cm}$&
$\hspace*{0.8cm}   L^{2}_\kappa(D)   \hspace*{0.8cm}$ &
$\hspace*{0.8cm}   H^{1}_\kappa(D)  \hspace*{0.8cm}$\\
\hline\hline
       $121$       &  $2.45$    & $27.25$ &  $2.38$    & $26.13$  \\
\hline
      $202$      & $0.63$    & $14.43$ &  $0.54$    & $12.29$\\
\hline
      $283$      & $0.36$    & $11.74$ &  $0.28$    & $9.01$\\
\hline

      $364$      & $0.29$    & $10.61$ &  $0.21$    & $7.49$\\
\hline
       $445$    & $0.16$  &$7.40$  &  $--$    & $--$\\
\hline
\end{tabular}
\caption{Convergence history for GMsFEM with different coarse spaces dimensions corresponding to the permeability field in Fig. \ref{fig:frac_various_exp2}.}
\label{table:frac_various_exp2}
\end{table}
\begin{table}[htb]
\centering
\begin{tabular}{|c|c|c|c|c|c|c|}
\hline
\multirow{2}{*}{$\text{dim}(V_{\text{off}})$} &
\multicolumn{2}{c|}{  $\|u-u_{\text{off}} \|$ (\%) } &
\multicolumn{2}{c|}{  $\|u_{\text{snap}}-u_{\text{off}} \|$ (\%) }\\
\cline{2-5} {}&
$\hspace*{0.8cm}   L^{2}_\kappa(D)   \hspace*{0.8cm}$ &
$\hspace*{0.8cm}   H^{1}_\kappa(D)  \hspace*{0.8cm}$&
$\hspace*{0.8cm}   L^{2}_\kappa(D)   \hspace*{0.8cm}$ &
$\hspace*{0.8cm}   H^{1}_\kappa(D)  \hspace*{0.8cm}$\\
\hline\hline
       $121$       &  $0.58$    & $11.77$ &  $0.46$    & $9.56$  \\
\hline
      $202$      & $0.28$    & $8.93$ &  $0.17$    & $5.74$\\
\hline

      $283$  & $0.24$  &$7.98$ &  $0.10$    & $4.11$\\
\hline
       $364$    & $0.17$  &$6.84$  &  $--$    & $--$\\
\hline
\end{tabular}
\caption{Convergence history for GMsFEM with different coarse spaces dimensions corresponding to the permeability field in Fig. \ref{fig:frac_various_exp3}.}
\label{table:frac_various_exp3}
\end{table}
%
\begin{table}[htb]
\centering
\begin{tabular}{|c|c|c|c|c|c|c|}
\hline
\multirow{2}{*}{$\text{dim}(V_{\text{off}})$} &
\multicolumn{2}{c|}{  $\|u-u_{\text{off}} \|$ (\%) } &
\multicolumn{2}{c|}{  $\|u_{\text{snap}}-u_{\text{off}} \|$ (\%) }\\
\cline{2-5} {}&
$\hspace*{0.8cm}   L^{2}_\kappa(D)   \hspace*{0.8cm}$ &
$\hspace*{0.8cm}   H^{1}_\kappa(D)  \hspace*{0.8cm}$&
$\hspace*{0.8cm}   L^{2}_\kappa(D)   \hspace*{0.8cm}$ &
$\hspace*{0.8cm}   H^{1}_\kappa(D)  \hspace*{0.8cm}$\\
\hline\hline
       $121$       &  $1.61$    & $24.46$ &  $1.53$    & $23.21$  \\
\hline
      $202$      & $0.45$    & $13.33$ &  $0.37$    & $10.98$\\
\hline

      $283$  & $0.27$  &$10.33$ &  $0.17$    & $7.07$\\
\hline

 $364$  & $0.20$  &$8.58$ &  $0.09$    & $4.13$\\
\hline
       $445$    & $0.17$  &$7.51$  &  $--$    & $--$\\
\hline
\end{tabular}
\caption{Convergence history for GMsFEM with different coarse spaces dimensions corresponding to the permeability field in Fig. \ref{fig:fracture_construction}.}
\label{table:fracture_construction}
\end{table}
%
\begin{table}[htb]
\centering
\begin{tabular}{|c|c|c|c|c|c|c|}
\hline
\multirow{2}{*}{$\text{dim}(V_{\text{off}})$} &
\multicolumn{2}{c|}{  $\|u-u_{\text{off}} \|$ (\%) } &
\multicolumn{2}{c|}{  $\|u_{\text{snap}}-u_{\text{off}} \|$ (\%) }\\
\cline{2-5} {}&
$\hspace*{0.8cm}   L^{2}_\kappa(D)   \hspace*{0.8cm}$ &
$\hspace*{0.8cm}   H^{1}_\kappa(D)  \hspace*{0.8cm}$&
$\hspace*{0.8cm}   L^{2}_\kappa(D)   \hspace*{0.8cm}$ &
$\hspace*{0.8cm}   H^{1}_\kappa(D)  \hspace*{0.8cm}$\\
\hline\hline
       $121$       &  $1.28$    & $19.54$ &  $1.22$    & $18.57$  \\
\hline
      $202$      & $0.25$    & $9.39$ &  $0.20$    & $7.23$\\
\hline

      $283$  & $0.19$  &$8.00$ &  $0.13$    & $5.31$\\
\hline

 $364$      & $0.16$    & $7.10$ &  $0.09$    & $3.82$\\
\hline
       $445$    & $0.12$  &$5.98$  &  $--$    & $--$\\
\hline
\end{tabular}
\caption{Convergence history for GMsFEM with different coarse spaces dimensions corresponding to the permeability field in Fig. \ref{fig:fracture_construction1}.}
\label{table:fracture_construction1}
\end{table}
We have also performed a few simulations with uniform triangular meshes with the fracture field depicted in Fig. \ref{fig:frac_tri}. We depict the fine-scale velocity components in Fig. \ref{fig:tri_velocity}. The coarse-scale velocity component of this fracture field is similar to those fracture fields above, and thus are not presented. Triangular meshes can also be used in the simulation.

\subsection{Numerical results with EFM}\label{subsec:efm}

\begin{figure}[htb]
\centering
 \subfigure[Fracture field]{\label{fig:frac_hierarchical}
    \includegraphics[width = 0.30\textwidth, keepaspectratio = true]{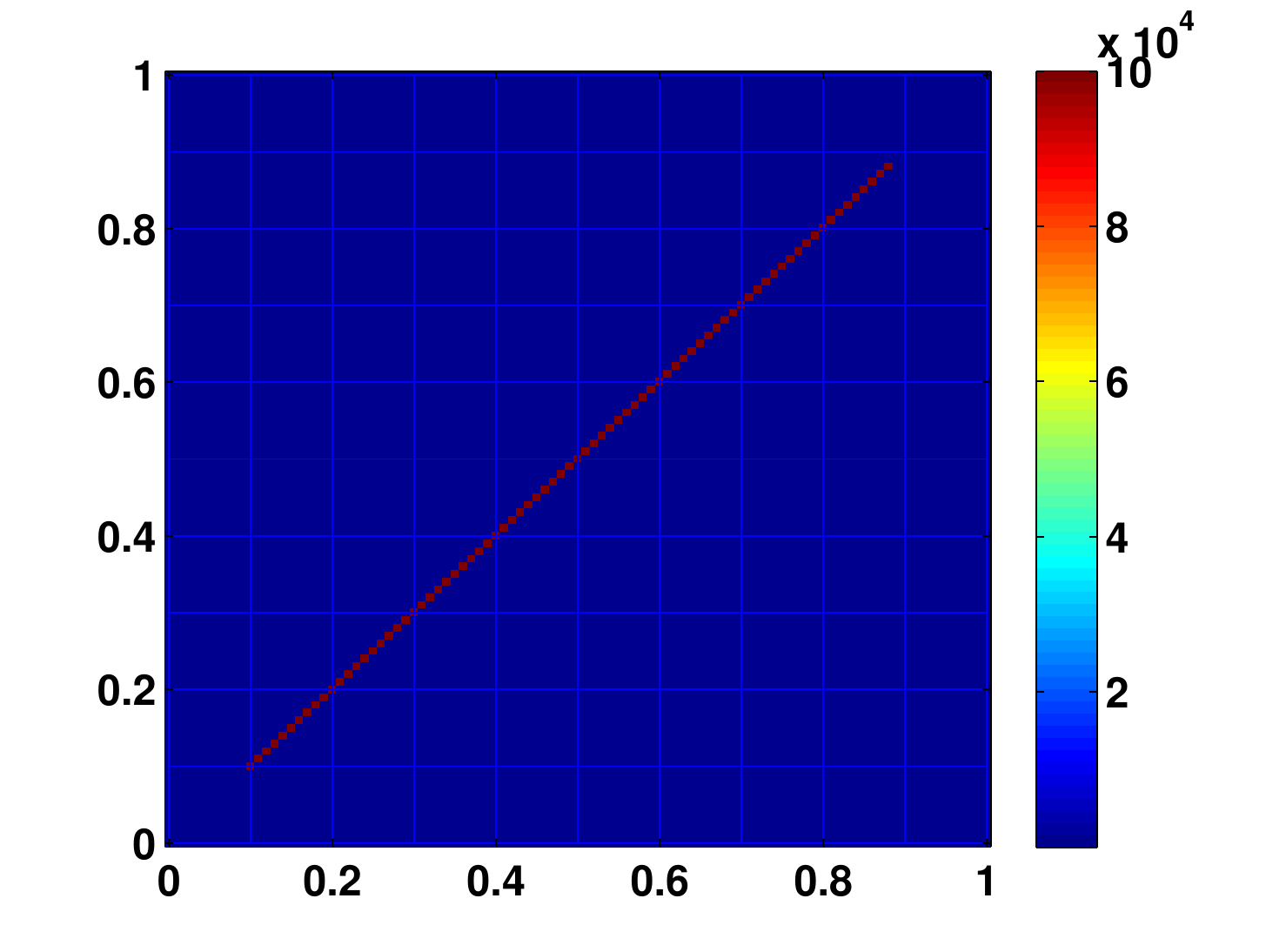}
   }
  \subfigure[Fine-scale solution]{\label{fig:FEMsoln_frac_hierarchical}
    \includegraphics[width = 0.30\textwidth, keepaspectratio = true]{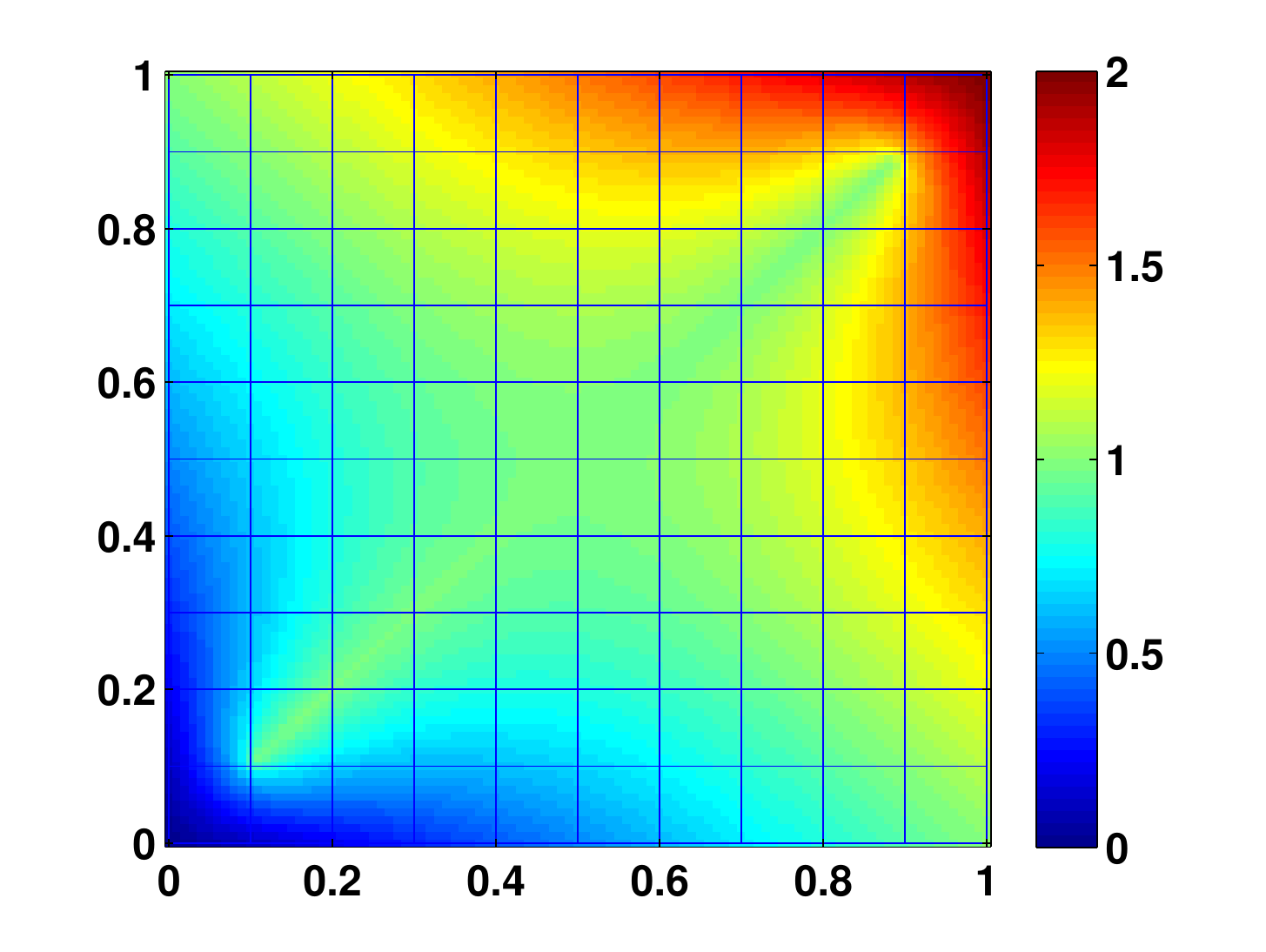}
   }
 \subfigure[coarse-scale solution]{\label{fig:umsfine_hierarchical}
    \includegraphics[width = 0.30\textwidth, keepaspectratio = true]{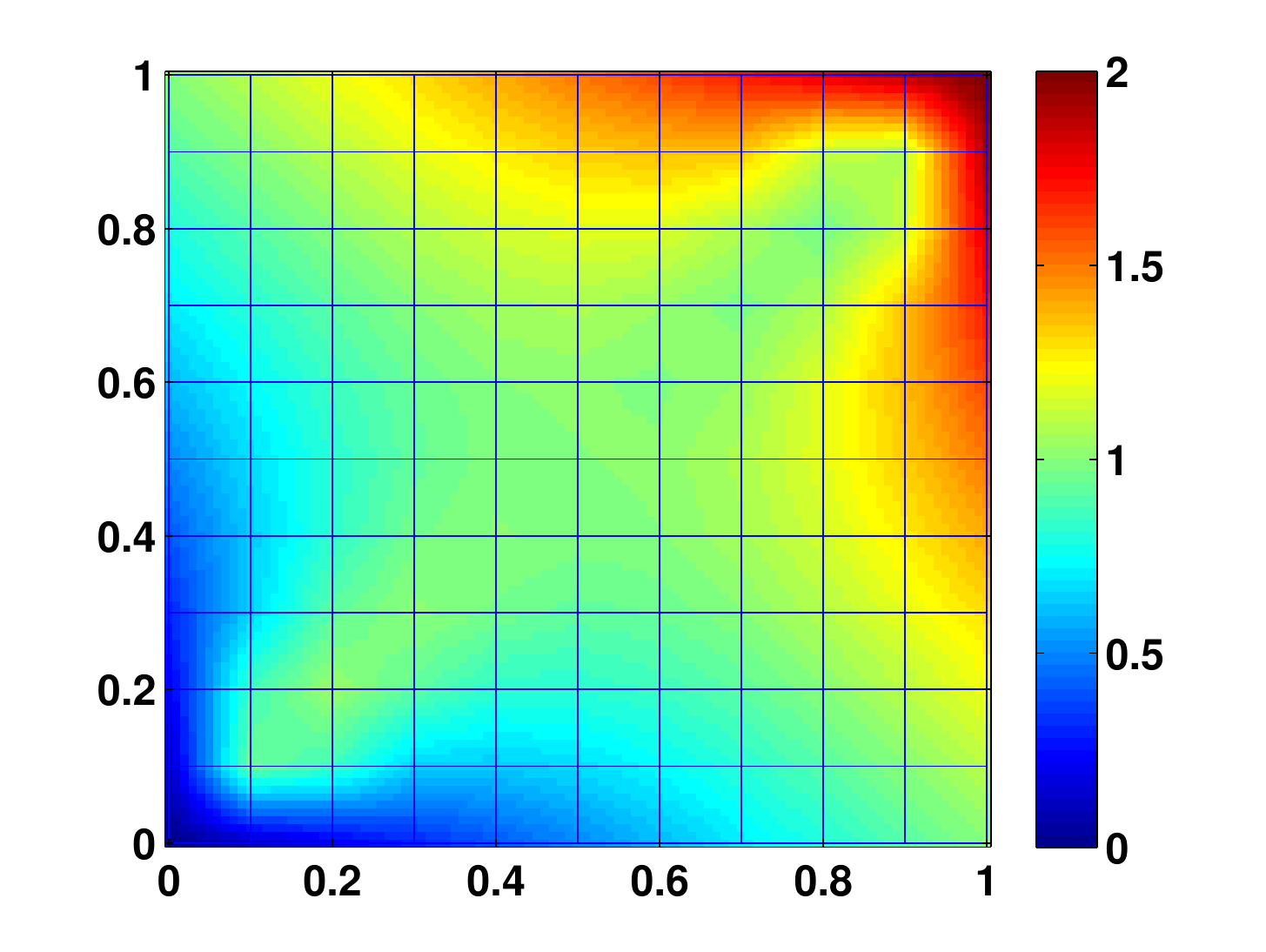}
   }
 \caption{Numerical results using GMsFEM with EFM.}
 \label{fig:test_hierarchical}
\end{figure}

\begin{figure}[htb]
\centering
 \subfigure[Fracture field]{\label{fig:frac_hierarchical2}
    \includegraphics[width = 0.3\textwidth, keepaspectratio = true]{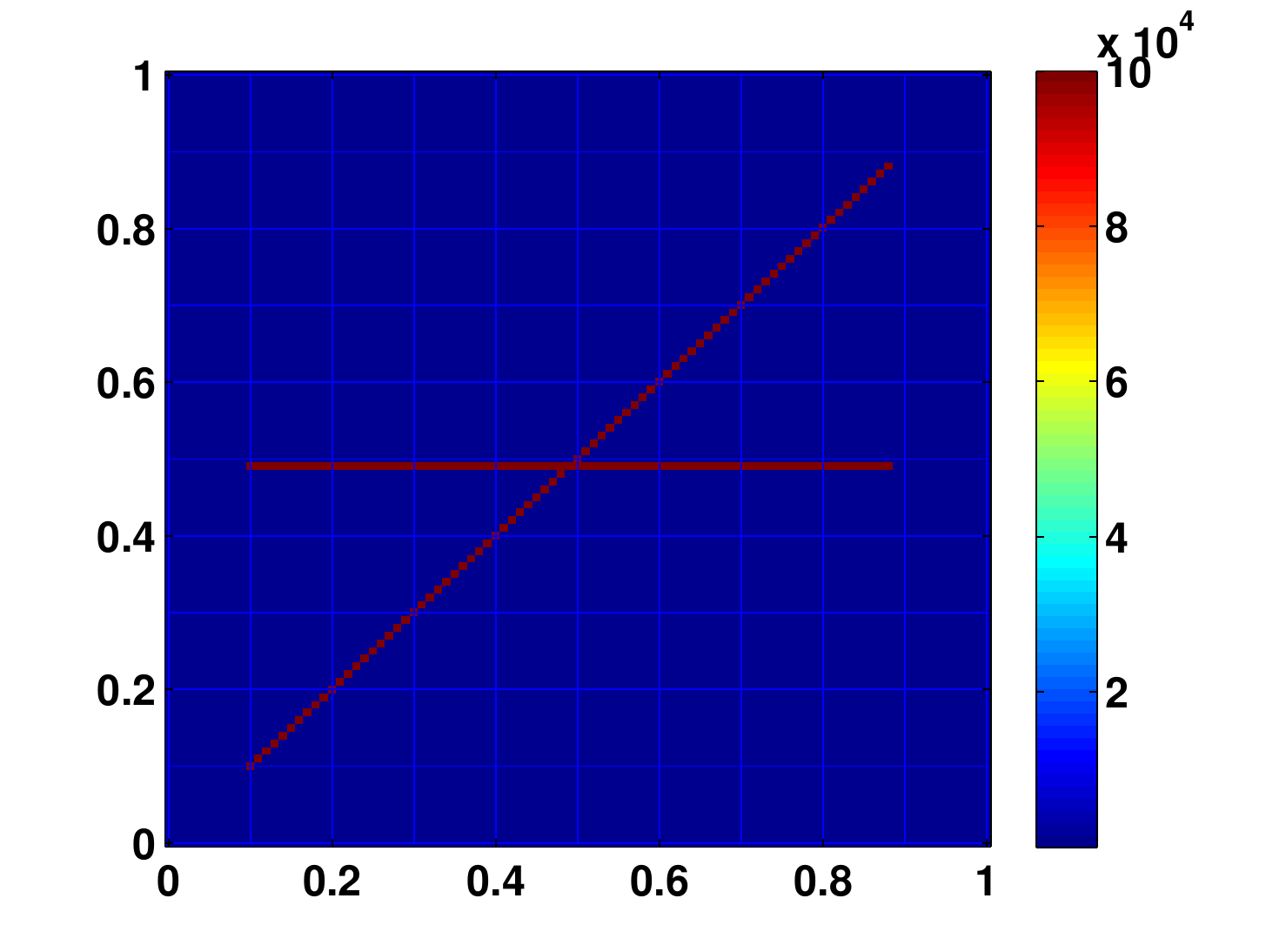}
   }
  \subfigure[Fine-scale solution]{\label{fig:FEMsoln_frac_hierarchical2}
    \includegraphics[width = 0.3\textwidth, keepaspectratio = true]{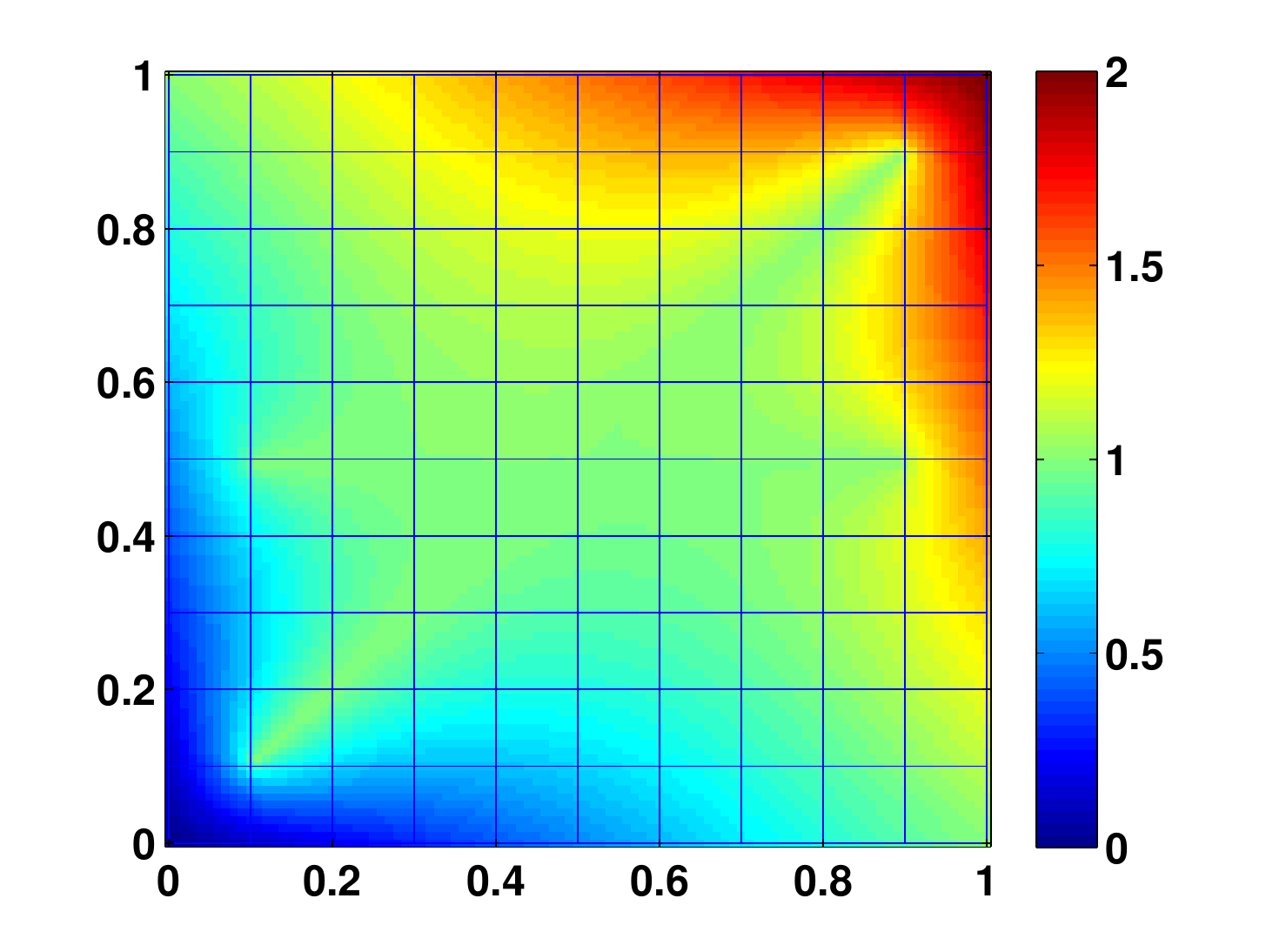}
   }
 \subfigure[coarse-scale solution]{\label{fig:umsfine_hierarchical2}
    \includegraphics[width = 0.3\textwidth, keepaspectratio = true]{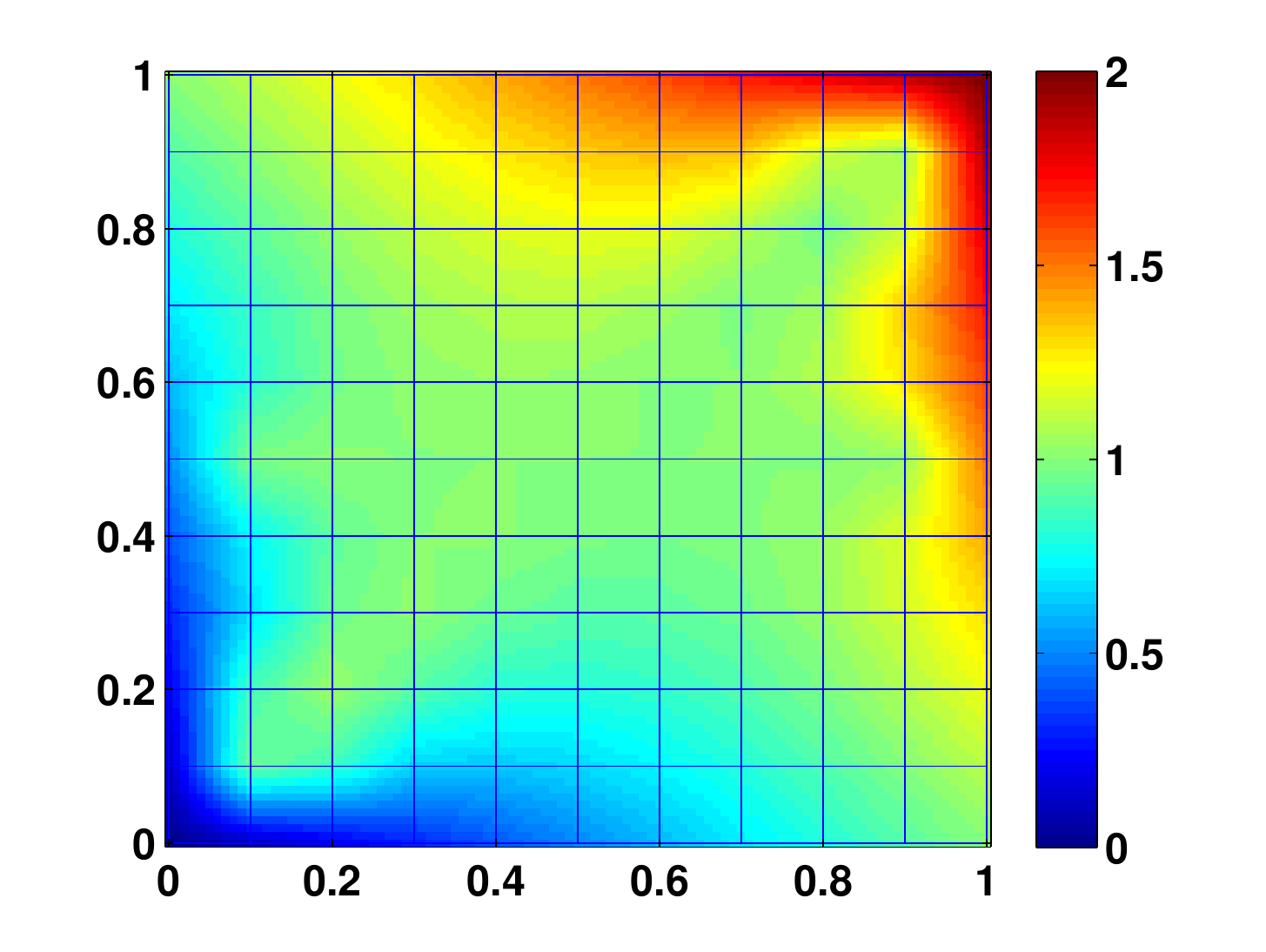}
   }
 \caption{Numerical results using GMsFEM with EMF.}
 \label{fig:test_hierarchical2}
\end{figure}

In this subsection, the numerical results using EFM are presented. As mentioned above, the main advantage of this model is that the fracture does not need to align with fine grid boundaries.
The numerical results are shown in Figs. \ref{fig:test_hierarchical} and \ref{fig:test_hierarchical2}.
Note that for this model, we have unique nodes (independent of the matrix nodes) for the
fractures. The number of long fractures are quite few (e.g. there is only one fracture in Fig. \ref{fig:test_hierarchical}). Therefore, it is not necessary for the upscaling over the long fractures.
The classical multiscale finite element method is applied in Fig. \ref{fig:test_hierarchical}
to obtain the coarse-scale solution. The corresponding weighted $H^1$ error and $L^2$ errors are $0.3327\%$ and $0.1734\%$. Therefore, no spectral problem is needed for this fracture field. Similarly,
for the simulation depicted in Fig. \ref{fig:test_hierarchical2}, the GMsFEM algorithm gives an error of $0.1344\%$ and $0.1409\%$ for energy error and $L^2$, error respectively.


\subsection{Numerical results for mutiscale DFM and EFM approach}
\label{subsec:dfm_efm}

\begin{figure}[htb]
\centering
 \subfigure[Fracture system]{\label{fig:frac_combined}
    \includegraphics[width = 0.30\textwidth, keepaspectratio = true]{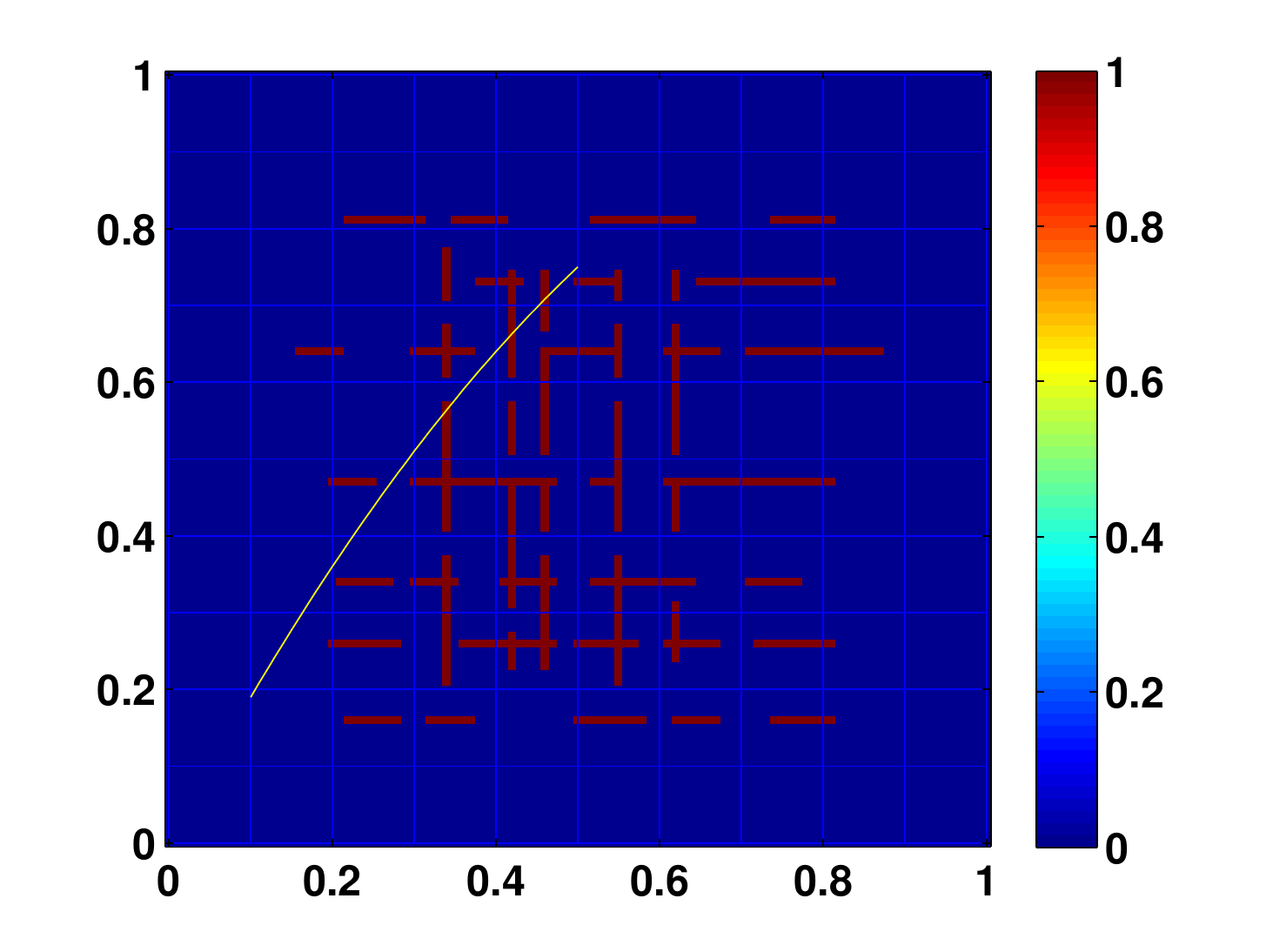}
   }
  \subfigure[Coarse-scale solution ]{\label{fig:umsfineSS_frac_combined}
    \includegraphics[width = 0.30\textwidth, keepaspectratio = true]{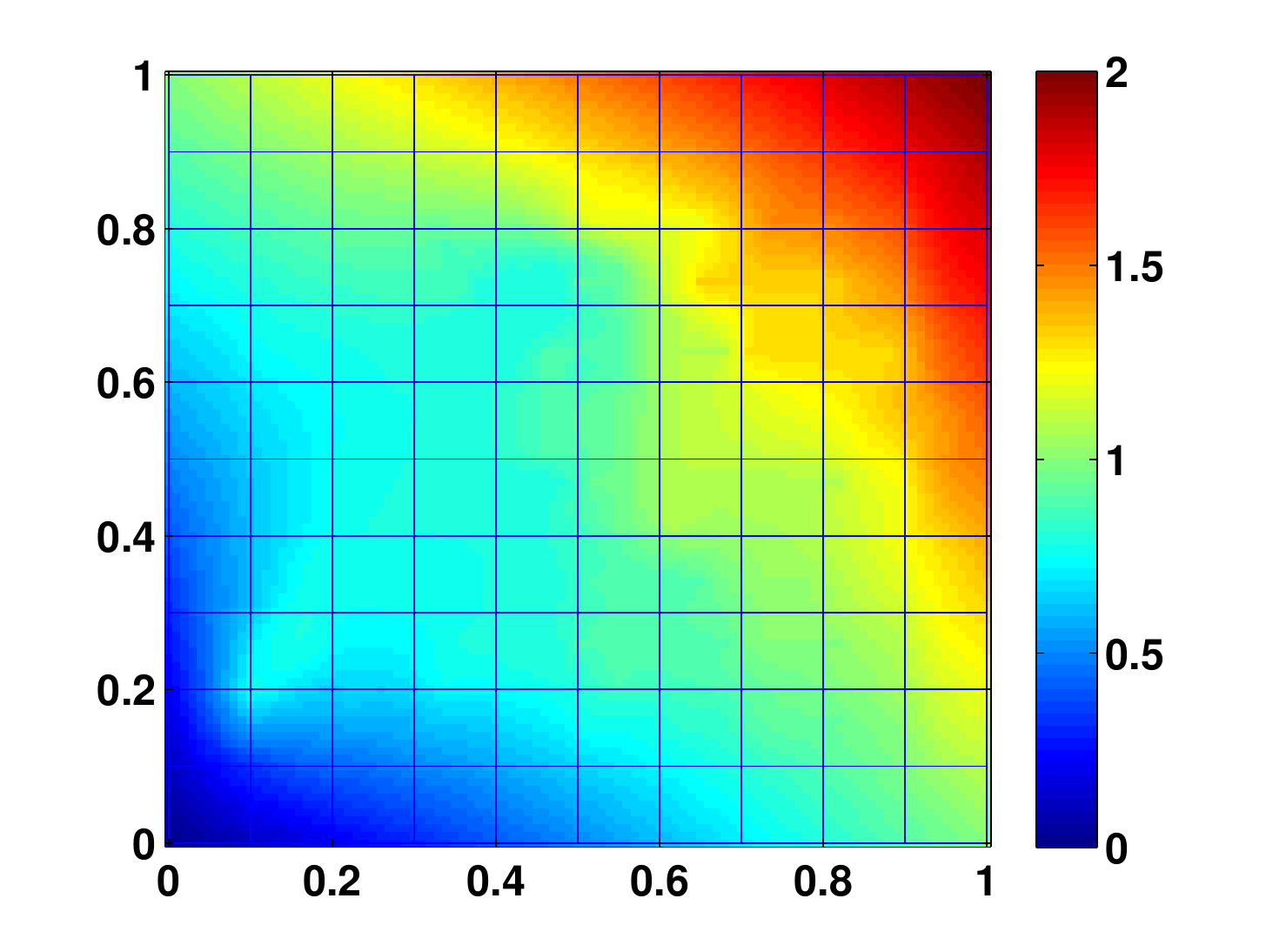}
   }
  \subfigure[Fine-scale solution]{\label{fig:FEMsoln_frac_combined}
    \includegraphics[width = 0.30\textwidth, keepaspectratio = true]{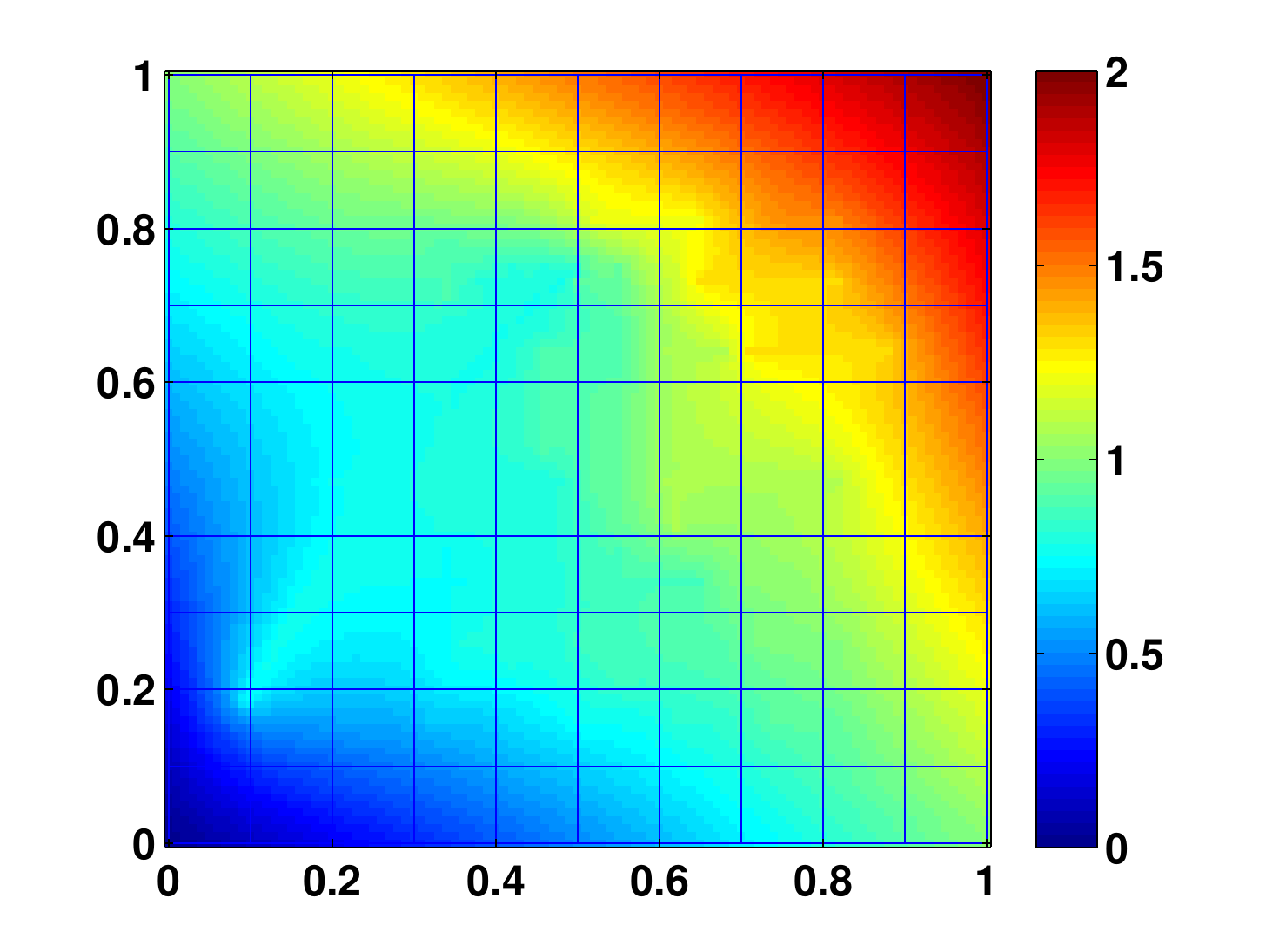}
   }
 \caption{Permeability field (left), coarse-scale solution (middle) and fine-scale solution (right) using the DFM model for short fractures and EFM for long fractures.}
 \label{fig:frac_combined_fine_coarse}
\end{figure}

In this subsection, we will
display several numerical results using DFM and EFM jointly.
 A heterogeneous medium with different lengths of fractures
is illustrated in Fig \ref{fig:frac_combined}. The
coarse-scale solution and the fine-scale solution are shown in Figs \ref{fig:umsfineSS_frac_combined}
and \ref{fig:FEMsoln_frac_combined}, respectively. The fine-scale solution is
calculated by the
DFM over the fine mesh. We observe that the coarse-scale solution is
a good approximation of
the fine-scale solution. Besides, the relative energy error and $L^2$
relative errors are
shown in Table \ref{table:fracture_combined}.  We see that the energy
relative error and
$L^2$ relative error are $9.17\%$ and $0.02\%$ when the dimension of
 the coarse system is 364.

%
\begin{table}[htb]
\centering
\begin{tabular}{|c|c|c|c|c|c|c|}
\hline
\multirow{2}{*}{$\text{dim}(V_{\text{off}})$} &
\multicolumn{2}{c|}{  $\|u-u_{\text{off}} \|$ (\%) } &
\multicolumn{2}{c|}{  $\|u_{\text{snap}}-u_{\text{off}} \|$ (\%) }\\
\cline{2-5} {}&
$\hspace*{0.8cm}   L^{2}_\kappa(D)   \hspace*{0.8cm}$ &
$\hspace*{0.8cm}   H^{1}_\kappa(D)  \hspace*{0.8cm}$&
$\hspace*{0.8cm}   L^{2}_\kappa(D)   \hspace*{0.8cm}$ &
$\hspace*{0.8cm}   H^{1}_\kappa(D)  \hspace*{0.8cm}$\\
\hline\hline
        $121$       &  $3.14$    & $96.16$ &  $3.05$    & $85.15$  \\
\hline
      $202$      & $2.22$    & $73.10$ &  $2.15$    & $63.04$\\
\hline
      $283$      & $0.04$    & $12.59$ &  $0.03$    & $5.66$\\
\hline

 $364$  & $0.02$  &$9.17$ &  $0.01$    & $2.41$\\
\hline
       $445$    & $0.01$  &$6.66$  &  $--$    & $--$\\
\hline
\hline
\end{tabular}
\caption{Convergence history for GMsFEM with different coarse spaces dimensions corresponding to the permeability field in Fig. \ref{fig:frac_combined_fine_coarse} using DFM for resolving the short fractures and hierarchical method for the long fractures.}
\label{table:fracture_combined}
\end{table}

\begin{figure}[htb]
\centering
 \subfigure[Fracture system]{\label{fig:frac_combined1}
    \includegraphics[width = 0.30\textwidth, keepaspectratio = true]{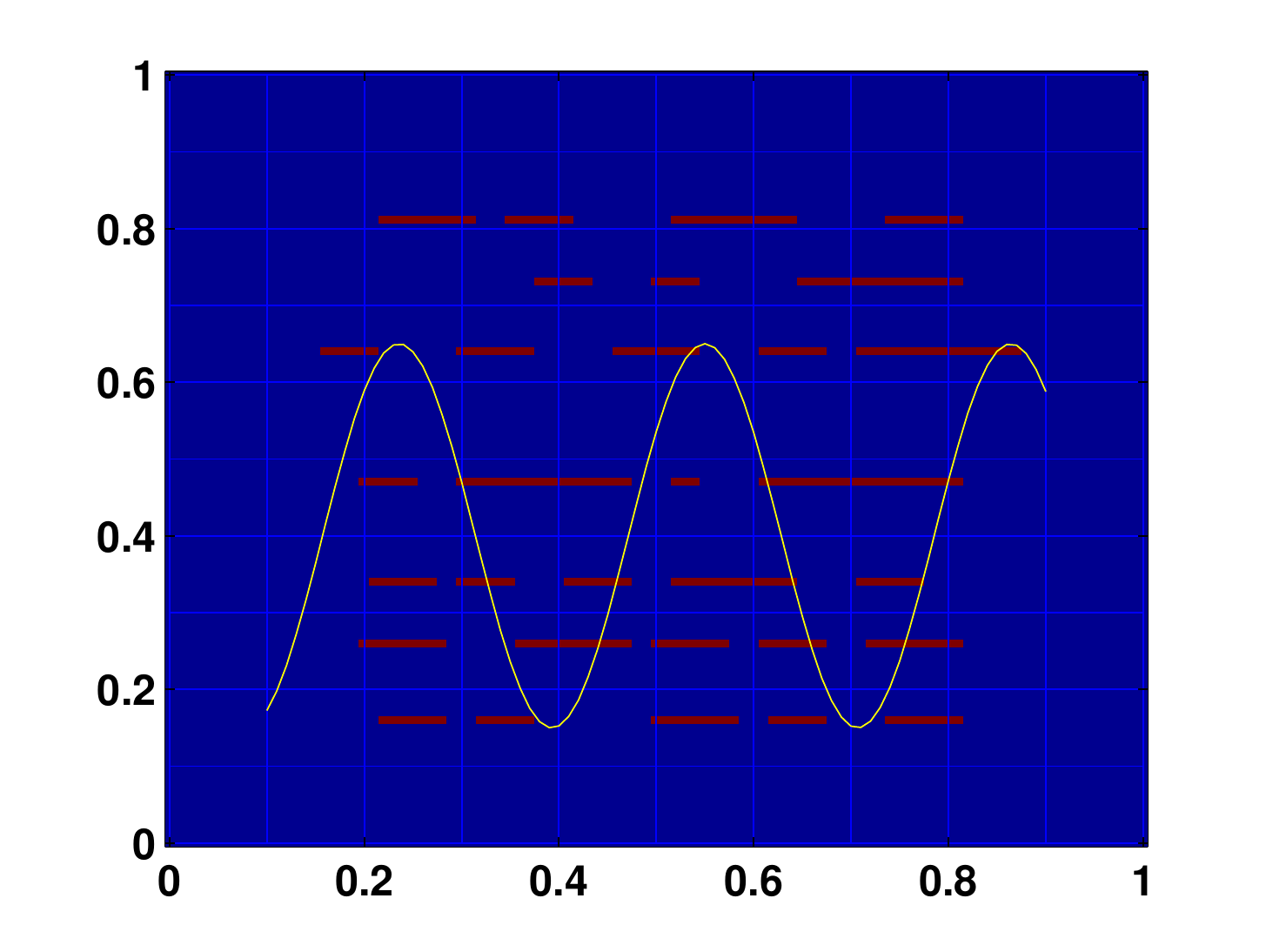}
   }
  \subfigure[Coarse-scale solution ]{\label{fig:umsfineSS_frac_combined1}
    \includegraphics[width = 0.30\textwidth, keepaspectratio = true]{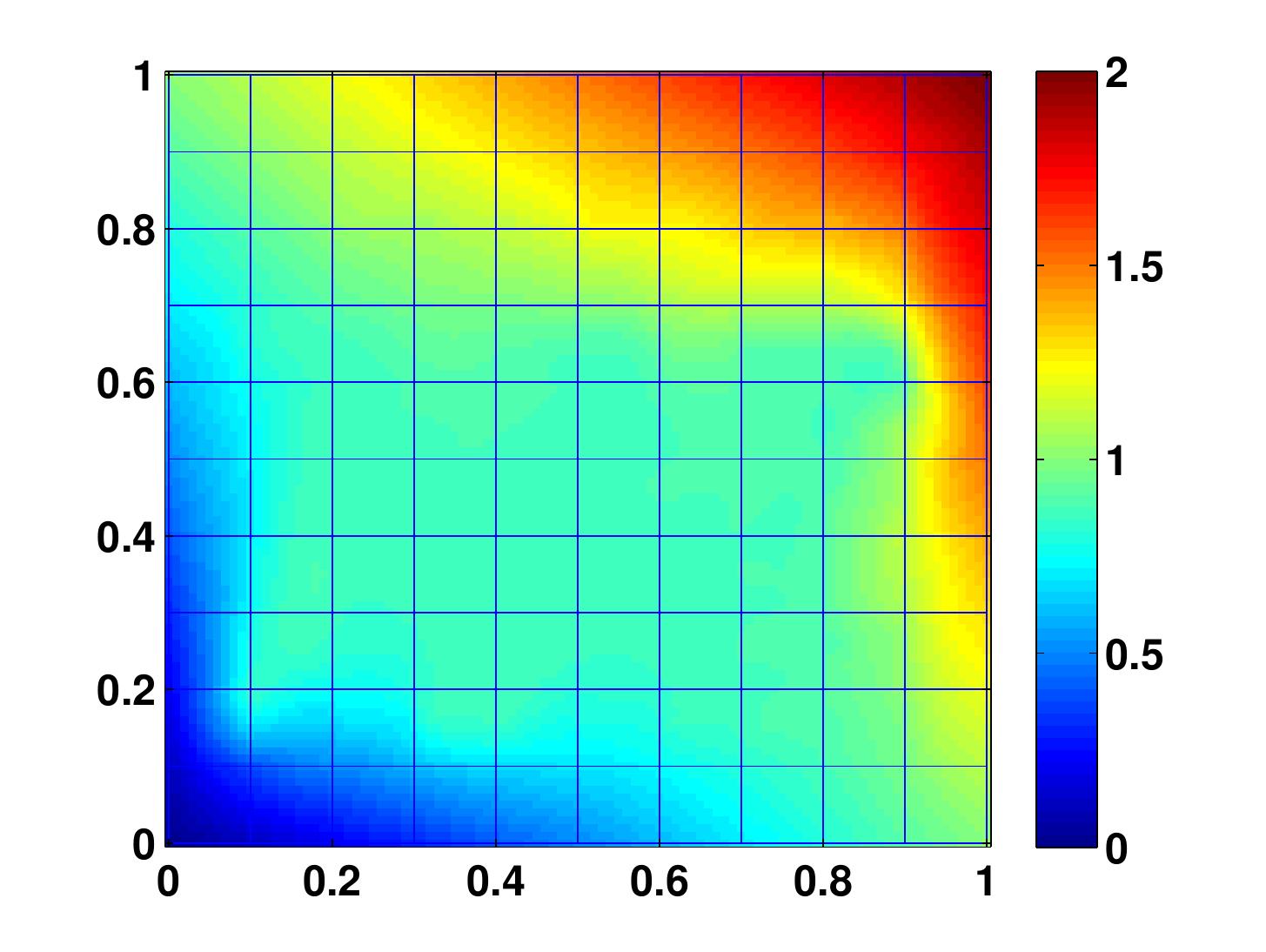}
   }
  \subfigure[Fine-scale solution]{\label{fig:FEMsoln_frac_combined1}
    \includegraphics[width = 0.30\textwidth, keepaspectratio = true]{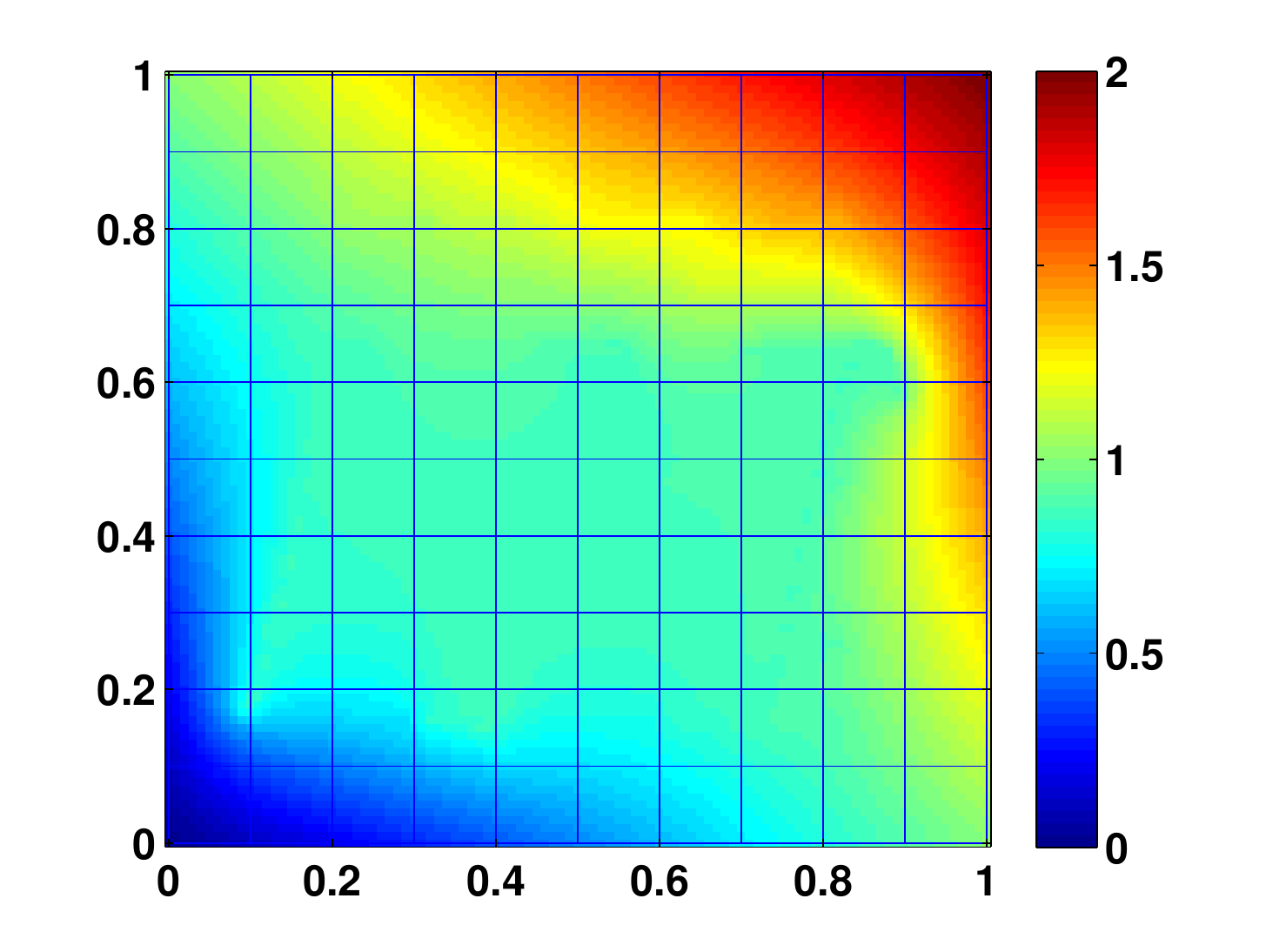}
   }
 \caption{Permeability field (left), coarse-scale solution (middle) and fine-scale solution (right) using the DFM model for short fractures and EFM for long fractures.}
 \label{fig:frac_combined1_fine_coarse}
\end{figure}

We test another heterogeneous fractured medium in Fig. \ref{fig:frac_combined1}.
Compared with the previous heterogeneous medium, this one has a longer fracture with a curved shape.
The coarse-scale solution and fine-scale solution are shown in Figs. \ref{fig:umsfineSS_frac_combined1}
and \ref{fig:FEMsoln_frac_combined1}.
We present the errors in Table \ref{table:fracture_combined1}.
\begin{table}[htb]
\centering
\begin{tabular}{|c|c|c|c|c|c|c|}
\hline
\multirow{2}{*}{$\text{dim}(V_{\text{off}})$} &
\multicolumn{2}{c|}{  $\|u-u_{\text{off}} \|$ (\%) } &
\multicolumn{2}{c|}{  $\|u_{\text{snap}}-u_{\text{off}} \|$ (\%) }\\
\cline{2-5} {}&
$\hspace*{0.8cm}   L^{2}_\kappa(D)   \hspace*{0.8cm}$ &
$\hspace*{0.8cm}   H^{1}_\kappa(D)  \hspace*{0.8cm}$&
$\hspace*{0.8cm}   L^{2}_\kappa(D)   \hspace*{0.8cm}$ &
$\hspace*{0.8cm}   H^{1}_\kappa(D)  \hspace*{0.8cm}$\\
\hline\hline
        $121$       &  $0.63$    & $37.85$ &  $0.54$    & $28.69$  \\
\hline
      $202$      & $0.52$    & $32.23$ &  $0.45$    & $22.55$\\
\hline
      $283$      & $0.03$    & $10.15$ &  $0.01$    & $2.54$\\
\hline

 $364$  & $0.02$  &$8.13$ &  $0.0008$    & $0.09$\\
\hline
       $445$    & $0.02$  &$7.55$  &  $--$    & $--$\\
\hline
\hline
\end{tabular}
\caption{Convergence history for GMsFEM with different coarse spaces dimensions corresponding to the permeability field in Fig. \ref{fig:frac_combined1_fine_coarse} using DFM for resolving the short fractures and EFM for long fractures.}
\label{table:fracture_combined1}
\end{table}


\subsection{Adaptive method}
\label{subsec:adaptive}

In this subsection, we will discuss adaptive strategies for generating
multiscale basis functions efficiently.
First, we will not use extra basis functions in the regions with long
fractures. This intuitive approach is based on previous findings
\cite{egw10}. Next,
we will use the error indicators that we proposed in \cite{Chung2013},
which identify the regions adaptive enrichment is needed.

We consider the fracture fields shown in Figs. \ref{fig:frac_various_exp2} and \ref{fig:frac_various_exp3}. As shown in these figures, both fields have channels and isolated inclusions in certain coarse neighborhoods. The local coarse spaces contain enough information about isolated inclusions because of the usage of multiscale partition of unity functions in the construction of global generalized multiscale basis. Therefore, the adaptive enrichment process takes place in the regions with long channels.

We begin with $1$ basis in each coarse node and examine the fracture field shown in Fig.\ref{fig:frac_various_exp2} following the reasoning above. We conclude that the coarse nodes with adaptive enrichment are $(i,j)$, $3\leq i\leq 9$ and $5\leq j\leq 9$. Adding 3 more basis to those coarse nodes, we can get the energy error of $8.55\%$ with the coarse space dimension of $226$. Compare with results listed in Table \ref{table:frac_various_exp2} (coarse space of dimension $364$ and the energy error of $10.61\%$), we observe that with this adaptive enrichment, we can get a better offline solution with a smaller solution space.

We also test the field  in Fig. \ref{fig:frac_various_exp3}. The coarse nodes with long fractures are $(i,j)$, $3\leq i\leq 9$ and $6\leq j\leq 8$. Adding 2 more basis to those coarse nodes, we can get the energy error of $7.70\%$ with the coarse space dimension of $163$. Compare with results listed in Table \ref{table:frac_various_exp3}, the energy error is $7.98\%$ with a coarse space of dimension $283$. We can see that with this adaptive algorithm, we can get a better offline solution with a smaller solution space.

Next, we use the $H^{-1}$-based error indicator proposed in \cite{Chung2013}.
We start with $1$ basis per coarse node and take $\theta=0.7$ in the adaptive algorithm proposed in \cite{Chung2013}. Comparing Table \ref{table:fracture_construction} with Table \ref{table:fracture_construction_adptv}, we notice that there is a slight gain regarding to the energy error
for this adaptive algorithm. In Table \ref{table:fracture_construction}, the energy error is $4.51\%$ for an offline space of dimension $445$, while $4.39\%$
for $424$ in the adaptive algorithm.

Next, we study the adaptive algorithm for the field in Fig. \ref{fig:frac_various_exp3}. From results in Table \ref{table:frac_various_exp3} and
Table \ref{table:fracture_construction1_adptv}, we conclude that the adaptive algorithm can improve the performance of GMsFEM. Using an offline space with $80$ less basis, the adaptive enrichment algorithm can still have a better performance.

%
\begin{table}[htb]
\centering
\begin{tabular}{|c|c|c|c|c|c|c|}
\hline
\multirow{2}{*}{$\text{dim}(V_{\text{off}})$} &
\multicolumn{2}{c|}{  $\|u-u_{\text{off}} \|$ (\%) } &
\multicolumn{2}{c|}{  $\|u_{\text{snap}}-u_{\text{off}} \|$ (\%) }\\
\cline{2-5} {}&
$\hspace*{0.8cm}   L^{2}_\kappa(D)   \hspace*{0.8cm}$ &
$\hspace*{0.8cm}   H^{1}_\kappa(D)  \hspace*{0.8cm}$&
$\hspace*{0.8cm}   L^{2}_\kappa(D)   \hspace*{0.8cm}$ &
$\hspace*{0.8cm}   H^{1}_\kappa(D)  \hspace*{0.8cm}$\\
\hline\hline
       $121$       &  $1.61$    & $24.46$ &  $1.53$    & $23.21$  \\
\hline
      $159$      & $0.44$    & $13.49$ &  $0.36$    & $11.15$\\
\hline

      $230$  & $0.27$  &$10.36$ &  $0.20$    & $7.32$\\
\hline

 $362$  & $0.18$  &$8.37$ &  $0.14$    & $5.50$\\
\hline
       $424$    & $0.15$  &$7.45$  &  $--$    & $--$\\
\hline
\end{tabular}
\caption{Convergence history for GMsFEM with different coarse spaces dimensions using adaptive algorithm  corresponding to permeability field in Fig. \ref{fig:fracture_construction}.}
\label{table:fracture_construction_adptv}
\end{table}
%
\begin{table}[htb]
\centering
\begin{tabular}{|c|c|c|c|c|c|c|}
\hline
\multirow{2}{*}{$\text{dim}(V_{\text{off}})$} &
\multicolumn{2}{c|}{  $\|u-u_{\text{off}} \|$ (\%) } &
\multicolumn{2}{c|}{  $\|u_{\text{snap}}-u_{\text{off}} \|$ (\%) }\\
\cline{2-5} {}&
$\hspace*{0.8cm}   L^{2}_\kappa(D)   \hspace*{0.8cm}$ &
$\hspace*{0.8cm}   H^{1}_\kappa(D)  \hspace*{0.8cm}$&
$\hspace*{0.8cm}   L^{2}_\kappa(D)   \hspace*{0.8cm}$ &
$\hspace*{0.8cm}   H^{1}_\kappa(D)  \hspace*{0.8cm}$\\
\hline\hline
       $121$       &  $1.28$    & $19.54$ &  $1.22$    & $18.57$  \\
\hline
      $169$      & $0.25$    & $9.35$ &  $0.20$    & $7.14$\\
\hline

      $259$  & $0.19$  &$8.01$ &  $0.13$    & $5.27$\\
\hline

 $310$      & $0.14$    & $6.85$ &  $0.09$    & $3.72$\\
\hline
       $370$    & $0.10$  &$5.49$  &  $--$    & $--$\\
\hline
\end{tabular}
\caption{Convergence history for GMsFEM with different coarse spaces dimensions using adaptive algorithm corresponding to permeability field in Fig. \ref{fig:fracture_construction1}.}
\label{table:fracture_construction1_adptv}
\end{table}

\subsection{Randomized snapshots}\label{subsec:randomized}

In this section,
we will investigate the oversampling randomized algorithm proposed in
\cite{randomized2014}. The main advantage of this algorithm is that
it uses much fewer snapshot basis functions to calculate the offline space.
Besides, the oversampling strategy is used to reduce the
subgrid errors due to
boundary conditions imposed to obtain the snapshot basis.
We refer to Fig. \ref{schematic} for an illustration of oversampling domain.

In our simulation, we set the oversampling size $t=2$ (two fine-grid layer
around the coarse region) and buffer number $p_{\text{bf}}^{\omega_i}=2$
for each coarse neighborhood $\omega_i$. I.e., we
 generate $n+2$ snapshot functions in each coarse neighborhood
when computing $n$ basis functions.
For the sake of completeness, we list the algorithm here
in Table {\ref{algorithm:random_snapshots}}.

 \begin{table}[h]
 \caption{Randomized GMsFEM Algorithm (\cite{randomized2014})}
 \begin{tabular}{r l}
 \hline\hline
 \\
    \textbf{Input}:& Fine grid size $h$, coarse grid size $H$,
oversampling size $t$, buffer number $p_{\text{bf}}^{\omega_i}$ for each $\omega_i$, \\
    & the number of local basis functions $k_{\text{nb}}^{\omega_i}$ for each $\omega_i$;\\
    \textbf{output}: & Coarse-scale solution $u_H$.\\

  1.& Generate oversampling region for each coarse block: $\mathcal{T}^H$, $\mathcal{T}^h$, and $\omega_i^{+}$; \\

  2.& Generate $k_{\text{nb}}^{\omega_i}+p_{\text{bf}}^{\omega_i}$ random vectors $r_l$ and obtain randomized snapshots in $\omega_i^{+}$;  \\

& Add a snapshot that represents the constant function on $\omega_i^+$;\\
  3. & Obtain $k_{\text{nb}}^{\omega_i}$ offline basis by a spectral decomposition restricted to the original domain $\omega_i$; \\

4. & Construct multiscale basis functions and solve it.\\

 \hline\hline
\label{algorithm:random_snapshots}
 \end{tabular}
 \end{table}
\begin{figure}[tb]
  \centering
  \includegraphics[width=1.0 \textwidth]{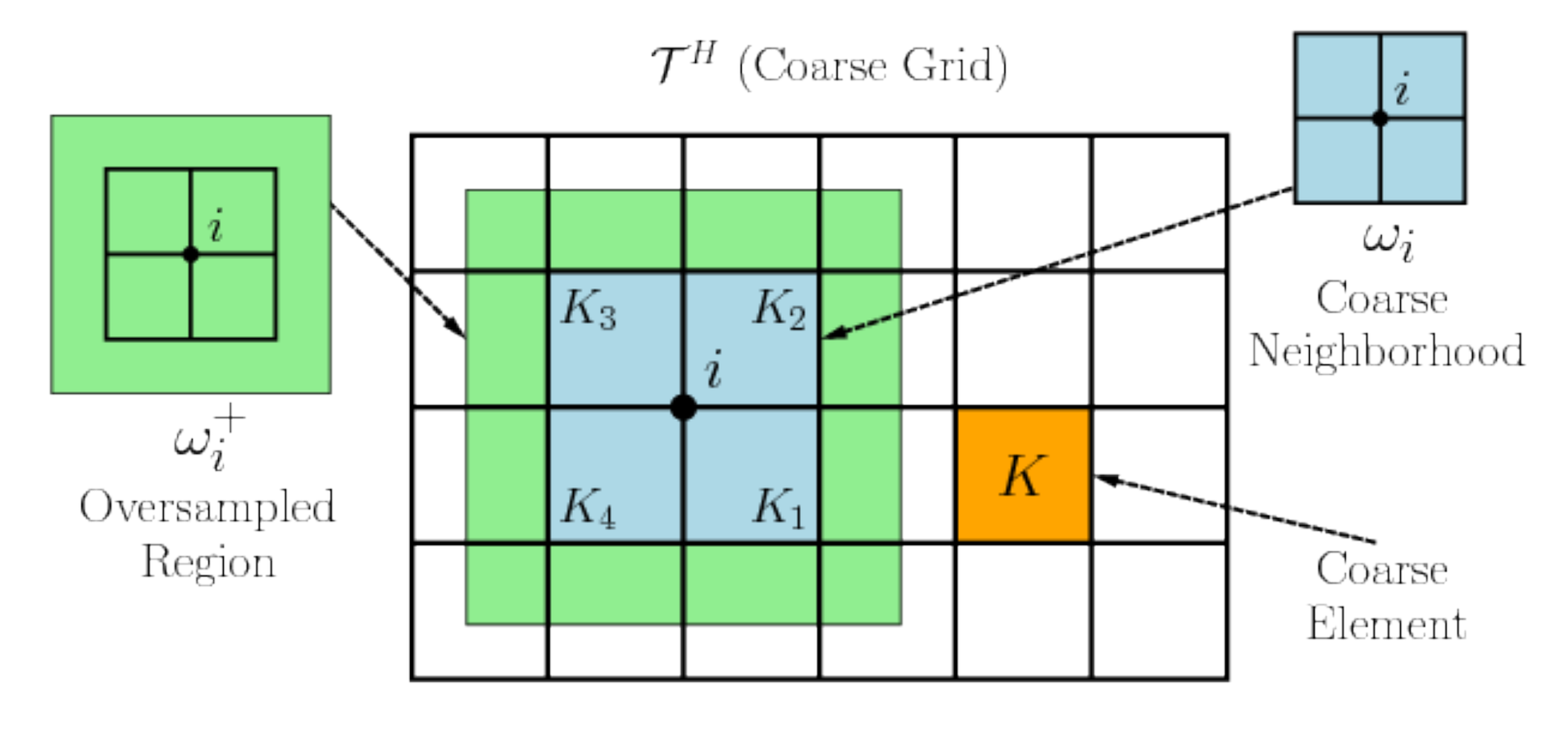}
  \caption{Illustration of a coarse neighborhood and oversampled domain. Here, $K$ is a coarse-grid block, $\omega_i$ is a coarse neighborhood of $x_i$, and
$\omega_i^+$ is an oversampled region}
  \label{schematic11}
\end{figure}

The numerical results using the permeability field
shown in Fig. \ref{fig:frac_combined1} are
presented in Table \ref{table:randomized}.
In this table, the first column shows the dimension
of the offline space, and the second columns displays the ratio between the
number of the randomized snapshot basis and the full snapshot basis. Then, in the following two columns, the weighted $L^2$ and energy error using the full snapshots are shown. At last, the results using the randomized snapshots alone is listed in the last two column.
Comparing these results with those using the full snapshot space, we deduce that the multiscale randomized algorithm performs well. The relative energy error is $15.82\%$ using the randomized snapshots, which account to $6.25\%$ of the full snapshots. However, we can only get $13.50\%$ relative energy error if the full snapshots are used instead.
\begin{table}[htb!]
\centering
\caption{Numerical results using randomized snapshots. $p_{\text{bf}}=4$.}
\label{table:randomized}
\small
\begin{tabular}{|r|c|c|c|c|c|c|c}
\hline
\multirow{2}{*}{$\text{dim}(V_{\text{off}})$}  &
\multirow{2}{*}{Snapshot ratio (\%)}  &
\multicolumn{2}{c|}{   All snapshots (\%) } &
\multicolumn{2}{c|}{   Randomized snapshots (\%) }\\
\cline{3-6} {}&{}&
$\hspace*{0.8cm}   L^{2}_\kappa(D)   \hspace*{0.8cm}$ &
$\hspace*{0.8cm}   H^{1}_\kappa(D)  \hspace*{0.8cm}$
&
$\hspace*{0.8cm}   L^{2}_\kappa(D)   \hspace*{0.8cm}$ &
$\hspace*{0.8cm}   H^{1}_\kappa(D)  \hspace*{0.8cm}$
\\
\hline\hline
      $283$ &$5.21$   & $7.31$ &    $15.55$    & $8.07$&    $16.85$    \\
\hline
      $364$ &$6.25$   &    $5.05$    & $13.50$& $6.60$ &    $15.82$    \\
\hline
       $445$&$7.29$   &   $4.46$    & $12.59$ &    $5.80$    & $14.48$\\
\hline
\end{tabular}
\end{table}

\section{Discussions and Conclusions}\label{sec:conclusion}

In this paper, we develop a multiscale finite element method for
flows in fractured
media. Our approach is based on Generalized
Multiscale Finite Element Method (GMsFEM).
Multiscale basis functions are constructed in the offline stage via local
spectral problems. To represent the fractures on the fine grid, we consider
two approaches (1)
Discrete Fracture Model (DFM) (2) Embedded Fracture Model (EFM)
and their combination.
The proposed procedure automatically selects multiscale basis functions
 via local spectral problems.

Numerical results are presented
where we demonstrate how one can adaptively add basis
functions in the regions of interest based on error indicators.
We consider various cases with long and short fractures.
Our numerical results show that GMsFEM with much fewer degrees of freedom compared to the
fine-scale simulations can capture the solution behavior accurately.
Furturemore.
we  discuss the use
of randomized snapshots (\cite{randomized2014}) which substantially reduces the offline computational cost.

\bibliographystyle{plain}
\bibliography{reference1}
\end{document}